\documentclass[10pt,reqno]{amsart} % please use amsart at 11pt

\usepackage{amssymb,latexsym}
\usepackage{cite} % to get Refs. [1,2,3] typeset as [1-3] automatically

\usepackage[height=190mm,width=130mm]{geometry} % this is the journal
\theoremstyle{plain}
\newtheorem{theorem}{Theorem}

\theoremstyle{definition}

\theoremstyle{remark}

\newcommand{\Z}{\mathbb{Z}}
\newcommand{\C}{\mathbb{C}}
\newcommand{\R}{\mathbb{R}}

\numberwithin{equation}{section} % to get equations numbered
\newcommand{\rmap}{\longrightarrow}

\newcommand{\U}{\ensuremath{\mathcal{U}}}

\newcommand{\F}{\ensuremath{\mathcal{F}}}

\newcommand{\g}{\ensuremath{\Gamma}}

\newcommand{\ps}{{\raise 1pt\hbox{\tiny (}}}

\newcommand{\pss}{{\raise 1pt\hbox{\tiny [}}}
\newcommand{\pdd}{{\raise 1pt\hbox{\tiny ]}}}
\newcommand{\pd}{{\raise 1pt\hbox{\tiny )}}}

\newcommand{\bs}{{\raise 1pt\hbox{\tiny [}}}
\newcommand{\bd}{{\raise 1pt\hbox{\tiny ]}}}

\def\cross{\mathinner{\mathrel{\raise0.8pt\hbox{$\scriptstyle>$}}
                 \joinrel\mathrel\triangleleft}}

\usepackage{citehack}
\usepackage{amsmath,amsthm}
\usepackage{amsfonts}
\usepackage{amssymb}
\usepackage{longtable}
\usepackage[matrix,arrow,curve]{xy}

\usepackage{lipsum}

\def\U{\mathcal{U}}

\def\K{\mathcal{K}}

\def\vH{{\check{H}}}

\usepackage{stackrel}

\newcommand{\be}{\begin{equation}}
\newcommand{\ee}{\end{equation}}

\newcommand{\nn}{\nonumber \\}

\newcommand{\wt}{\mbox{\rm wt}\ }

\newcommand{\nc}{\newcommand}
\nc{\cali}{\mathcal}
\nc{\on}{\operatorname}
\nc{\Wick}{{\mb :}}

\nc{\ddz}{\frac{\partial}{\partial z}}
\nc{\ch}{\mbox{ch}}
\nc{\Oo}{{\cali O}}
\nc{\cond}{|\,}
\nc{\bib}{\bibitem}
\nc{\pone}{\Pro^1}
\nc{\pa}{\partial}
\nc{\arr}{\rightarrow}
\nc{\larr}{\longrightarrow}
\nc{\ket}{\rangle}
\nc{\bra}{\langle}
\nc{\gam}{\bar{\gamma}}
\nc{\q}{\widetilde{Q}}
\nc{\ep}{\epsilon}
\nc{\su}{\widehat{{\mf s}{\mf l}}_2}
\nc{\sw}{{\mf s}{\mf l}}
\nc{\h}{{\mf h}}
\nc{\n}{{\mf n}}
\nc{\ab}{\mf{a}}
\nc{\is}{{\mb i}}
\nc{\js}{{\mb j}}
\nc{\bi}{\bibitem}
\nc{\He}{{\cali H}}
\nc{\inv}{^{-1}}
\nc{\ol}{\overline}
\nc{\wh}{\widehat}
\nc{\dst}{\displaystyle}

\nc{\delt}{\partial_t}
\nc{\ddt}{\frac{\partial}{\partial t}}
\nc{\delx}{\partial_x}
\nc{\mb}{\mathbf}
\nc{\mf}{\mathfrak}

\nc{\mbb}{\mathbb}
\nc{\Ctt}{\C((t))}
\nc{\Ct}{\C[t,t\inv]}

\nc{\ghat}{\wh{\g}}

\nc{\un}{\underline}
\nc{\mc}{\mathcal}
\nc{\BB}{{\mc B}}
\nc{\bb}{{\mf b}}
\nc{\kk}{{\mf k}}
\nc{\frob}{\times}
\nc{\sm}{\setminus}
\nc{\Pp}{{\mathbb P}^1}
\nc{\Aa}{{\mc A}}

\nc{\AutO}{\on{Aut}\Oo}
\nc{\AUTO}{\un{\on{Aut}}\Oo}
\nc{\AUTK}{\un{\on{Aut}}\K}
\nc{\Heout}{\He_{\out}}
\nc{\Hetil}{{\widetilde\He}}
\nc{\wb}{\overline}

\nc{\Res}{\on{Res}}
\nc{\pitil}{\Pi}
\nc{\Ctil}{\wt{C}}
\nc{\auto}{\on{Aut} \Oo}
\nc{\phitil}{\wt{\phi}}
\nc{\gz}{\g_{\vec z}}
\nc{\tensorM}{\bigotimes_{i=1}^N{\mathbb M}_i}
\nc{\tensorW}{\bigotimes_{i=1}^N W_{\nu_i,k}}
\nc{\out}{\on{out}}

\nc{\m}{{\mathfrak m}}

\nc{\gx}{\g^0_{\vec x}}

\nc{\hx}{\He^0_{\vec x}}
\nc{\tensorpi}{\pi_{\nu_1,\ldots,\nu_N}^\kappa}
\nc{\Phizw}{\Phi_{\vec w}({\vec z})}
\nc{\Pro}{{\mathbb P}}

\nc{\De}{\Delta}

\nc{\us}{\underset}
\nc{\Ll}{\mc L}
\nc{\dR}{\on{dR}}

\nc{\T}{{\mc T}}

\nc{\Xn}{\overset{\circ}X{}^n} \nc{\Dn}{\overset{\circ}D{}^n}
\nc{\Dxn}{\overset{\circ}D{}^n_x} \nc{\varphitil}{\wt{\varphi}}

\nc{\lf}{{\mf l}}
\nc{\GL}{{}^L G}
\nc{\Vir}{\on{Vir}}

\begin{document}
\title[Continual Lie algebras determined by chain complexes]  
{Continual Lie algebras determined by chain complexes} 
%%
 % please provide
                                % an abbreviated title 
%%\newcommand{\half}{\frac{1}{2}}
%%%%%%%%%%%%%%%%%%%%%%%%%%%%%%%%%%%%%%%%%%%%%%%%%%%%%%%%%%%%%%%%%%%%%%%%%%%%%%%
\author{A. Zuevsky} 
\address{Institute of Mathematics \\ Czech Academy of Sciences\\ Zitna 25, Prague\\Czech Republic}

\email{zuevsky@yahoo.com}

%%
%%%%%%%%%%%%%%%%%%%%%%%%%%%%%%%%%%%%%%%%%%%%%%%%%%%%%%%%%%%%%%%%%%%%%%%%%
% You may repeat \author \address as often as necessary                 %
%%%%%%%%%%%%%%%%%%%%%%%%%%%%%%%%%%%%%%%%%%%%%%%%%%%%%%%%%%%%%%%%%%%%%%%%%
%% 
\begin{abstract}
Continual Lie algebras
are infinite-dimensional generalizations of Lie algebras with discrete
root system by considering continual root systems. 
In this paper we establish the general relation 
between chain complexes and continual Lie algebras.  
%%  
%%%%%%%%%%%%%%%%%%%%%%%%%%%%%%%%%%%%%%%%%%%%%%%%%%%%%%%%%%%%%%%%%%%%%
%% 
The natural orthogonality condition with respect 
to a product among elements of a chain complex $\mathcal C$ spaces 
brings about to $\mathcal C$ the structure of a 
graded algebra with differential relations. 
%%   
%%%%%%%%%%%%%%%%%%%%%%%%%%%%%%%%%%%%%%%%%%%%%%%%%%%%%%%%%%%%%%%%%%%%
%%   
We prove the main result of this paper: a chain complex endowed 
with an appropriate Leibniz-property product of elements of its spaces 
and the Jacobi identity 
 brings about the structure of a 
continual Lie algebra
with the root space determined by parameters for the complex. 
That provides a new source of examples of continual Lie algebras. 
%% 
%%%%%%%%%%%%%%%%%%%%%%%%%%%%%%%%%%%%%%%%%%%%%%%%%%%%%%%%%%%%%%%%%%%%%
%% 
Finally, as an example, we consider the case of  
\v{C}ech-de Rham complex associated to a 
foliation of a smooth manifold. 
In a particular case of this chain complex, 
we derive explicitly the commutation relations 
for the corresponding continual Lie algebra. 
AMS Classification: 53C12, 57R20, 17B69 
\end{abstract}

\keywords{Continual Lie algebras, graded differential algebras, cohomology} 
\vskip12pt  % insert '\vskip12pt' while using '\twocolumn' command
%\vskip28pt % if there is no keywords

\maketitle

%%%%%%%%%%%%%%%%%%%%%%%%%%%%%%%%%%%%%%%%%%%%%%%%%%%%%%%%%%%%%%%%%%%%%%%%%%%%%%%%%%%%%%%%%%%%%%
%%%%%%%%%%%%%%%%%%%%%%%%%%%%%%%%%%%%%%%%%%%%%%%%%%%%%%%%%%%%%%%%%%%%%%%%%%%%%%%%%%%%%%%%%%%%%%
\begin{center}
{Conflict of Interest and Data availability Statements:}
\end{center}

The author states that: 

1.) The paper does not contain any potential conflicts of interests. 

2.) The paper does not use any datasets. No dataset were generated during and/or analysed 
during the current study. 

3.) The paper includes all data generated or analysed during this study. 

4.) Data sharing is not applicable to this article as no datasets were generated or analysed during the current study.

5.) The data of the paper can be shared openly.  

%%%%%%%%%%%%%%%%%%%%%%%%%%%%%%%%%%%%%%%%%%%%%%%%%%%%%%%%%%%%%%%%%%%%%%%%%%%%%%%%%%%%%%%%%%%%%%
%%%%%%%%%%%%%%%%%%%%%%%%%%%%%%%%%%%%%%%%%%%%%%%%%%%%%%%%%%%%%%%%%%%%%%%%%%%%%%%%%%%%%%%%%%%%%%%%%
%% 
\section{Introduction} 
\label{valued}
Continual Lie algebras introduced in \cite{save}
 are generalizations of
 infinite-dimensional Lie algebras with systems of discrete roots \cite{kac}.     
Then the notion was generalized and developed in \cite{sv1}-\cite{vershik}.  
These algebras are formulated in such 
a way that the space of roots is defined by continual 
sets of vectors.  
In commutation relations, 
the generators depend on kernels which are functionals of continual roots. 
 The Jacobi identity for continual Lie algebras 
result in non-trivial functional relations for kernels.   
Similar constructions appeared recently 
in the study of certain Hall algebras of coherent
sheaves described in terms of certain continuum limits of Kac-Moody
algebras \cite{ass}. 
%%

%%%%%%%%%%%%%%%%%%%%%%%%%%%%%%%%%%%%%%%%%%%%%%%%%%%%%%%%%%%%%%%%%
%%
Though the general theory of continual Lie algebras and their representations is missing, 
there exist many applications, especially in the theory of completely 
integrable and exactly solvable models 
\cite{razsav}.   
 Various applications can be also found in other fields of mathematics 
\cite{bakas, DS, vershik}.   
%%

%%%%%%%%%%%%%%%%%%%%%%%%%%%%%%%%%%%%%%%%%%%%%%%%%%%%%%%%%%%%%%%%%%%%%%%%%
%% 
In the original papers \cite{save}-\cite{vershik}, 
in order to find new classes of examples of continual 
Lie algebras
various approaches were studied. 
In a series of papers \cite{sv2, sv3} the authors have generalized first known examples and have found new non-trivial and 
fundamental ones arising from  
various branches of mathematics, in particular, functional analysis \cite{vershik}, 
differential geometry, 
algebraic geometry, non-commutative geometry,  
and mathematical physics \cite{bakas, FMc}.  
 Some of them came from physical applications,   
e.g., equations related to the Ricci flow and cosmology \cite{bakas, FMc}.   
Though they appear in various corners of modern mathematical research, at some stage, 
the arsenal of ideas of generation of  
new examples of continual Lie algebras was exhausted.  
An intriguing question is how to construct the general way  
which would allow to create continual Lie algebras.   
In many cases, this problem reduces to the question of finding 
new appropriate solutions for relations \eqref{jacjac} 
(which is a very interesting problem by itself)  
following from the Jacobi identity.  
%%

%%%%%%%%%%%%%%%%%%%%%%%%%%%%%%%%%%%%%%%%%%%%%%%%%%%%%%%%%%%%%%%%%%%%%%%%%%%%%%%%%%%%%%%%%%%%%%
%% 
The idea of this work is to use properties of the general chain complexes  
(with complex spaces depending on sets of parameters) endowed with a product 
in order to derive the structure of continual Lie algebras. 
%%%
The main result of this paper consists in the following  
%%  
%%%%%%%%%%%%%%%%%%%%%%%%%%%%%%%%%%%%%%%%%%%%%%%%%%%%%%%%%%%%%%%%%%%%%%%%%%%%%%%%%
%%
\begin{theorem}
\label{mainpro}
For the set $I_0$ of all pairs of independent elements 
 $\alpha_i^{(K_i)}$, 
  the orthogonality condition \eqref{ortho},   
with existence of the left and right principal ideals with 
respect to subset of the complex, 
 the generators 
$\left\{ \chi, \; \delta(n_0) \chi, \; \Phi, \; \delta(n) \Phi, \;   
\alpha_i^{(K_i)}, \;  
\delta(n(i))
 \alpha_i^{(K_i)} \right\}$,  
$n(i)= n_i^{(K_i)}$, 
 the relations \eqref{ogromno},   %\eqref{ser1}-\eqref{ser2} 
form a graded differential algebra. 
%%
%%%%%%%%%%%%%%%%%%%%%%%%%%%%%%%%%%%%%%%%%%%%%%%%%%%%%%%%%%%%%%%%%%%%%%%
%%
If, depending on properties of the chain complex 
$\left(C(n, \Theta_n), \cdot\right)$, 
 the kernels/mappings defining the commutation relations 
originating from the system \eqref{ogromno}  
satisfy the Jacobi identity \eqref{jacjac}, 
then the bigraded differential algebra turns into  
a continual Lie algebra $\mathcal G(\Theta_{n(i)})$ with the root space depending on 
the set of parameters $\Theta_{n(i)}$, of spaces $C\left(n(i), \Theta_{n(i)}\right)$.     
\end{theorem}
%%
%%%%%%%%%%%%%%%%%%%%%%%%%%%%%%%%%%%%%%%%%%%%%%%%%%%%%%%%%%%%%%%55
%%
The plan of the paper is the following. 
%%%
In the next subsection we recall the notion of a continual Lie algebra \cite{save}.  
In Section \ref{pseudo}, starting from the general (infinite) chain complex, 
we first define a product 
acting on elements of various spaces of a complex. 
%% 
%%%%%%%%%%%%%%%%%%%%%%%%%%%%%%%%%%%%%%%%%%%%%%%%%%%%%%%%%%%
%% 
We require the Leibniz rule to be fulfilled for 
(co)boundary operators with respect to this product.   
In addition to that we assume its skew-symmetry properties.  
Thus we endow the complex with the structure of a graded differential 
 algebra $\mathcal C$.  
%%

%%%%%%%%%%%%%%%%%%%%%%%%%%%%%%%%%%%%%%%%%%%%%%%%%%%%%%%%%%%%%%%%%%
%%
We then apply the natural condition of orthogonality 
with respect to a product among the chain complex spaces. 
For instance, in differential geometry, 
the orthogonality condition reduces to the integrability 
condition for differential forms. 
Further actions of (co)boundary operators,
as well as consequences from the definition of the 
multiplication, 
 generate a system of relations for elements of various chain complex spaces  
with extra compatibility conditions for indices.     
This determines the structure of a graded algebra with relations for $\mathcal C$.  
%%
%%%%%%%%%%%%%%%%%%%%%%%%%%%%%%%%%%%%%%%%%%%%%%%%%%%%%%%%%%%%%%%%%%%%5
%% 
By choosing independent elements satisfying the above relations,
 we then show the main result of this paper: 
for a graded differential algebra associated to a chain complex, 
the independent generators, the system of relations, 
and the Jacobi identity define the structure of  
a continual Lie algebra with 
the root space provided by the spaces parameters of a chain complex. 
%%% 

%%%%%%%%%%%%%%%%%%%%%%%%%%%%%%%%%%%%%%%%%%%%%%%%%%%%%%%%%%%%%%%%%%%%%%%%%%%%%%%%%%%%%%%%%%%%%%%%%%%%%%%%%%%%
%%
In Section \ref{example} we consider the example of
 the \v{C}ech-de Rham cochain bicomplex \cite{CM} (c.f. subsection \ref{cm})   
associated to a smooth manifold foliation \cite{DNF, BG, BGG, BEG}.    
This bicomplex has a deep geometric meaning \cite{galaev} and is defined for the spaces of 
differential forms and  
holonomy embeddings acting between sections of the transversal basis for a foliation.  
In subsection \ref{cm} we recall the notions of holonomy mappings, 
transversal sections, and transversal  
basis for a foliation.  
The space of differential forms has already a structure of bigraded differential algebra with respect to 
the natural product \cite{CM}.     
For the  bicomplex under consideration, we define the second product, 
 satisfying properties required for the construction of a system of  relations. 
 According to the exposition we give in Section \ref{pseudo}, 
 by involving the orthogonality conditions for an arbitrary pairs of the \v{C}ech-de Rham bicomplex, 
we apply the general scheme of construction of a bigraded differential algebra. 
Then we single out generators and commutation relations for the corresponding 
continual Lie algebra with the space of roots and kernels given by the sets of holonomy 
embeddings.    
In a particular geometric case associated to a codimension 
one foliation over a three-dimensional smooth manifold,  
we start from the integrability condition applied to elements of the same bicomplex space  
which leads to the appearance of 
a continual Lie algebra with generators, kernels, 
and commutation relations explicitly described. 
%%  
%%%%%%%%%%%%%%%%%%%%%%%%%%%%%%%%%%%%%%%%%%%%%%%%%%%%%%%%%%%%%%%%%%%%%%%%%%%%%%%%%
%%%%%%%%%%%%%%%%%%%%%%%%%%%%%%%%%%%%%%%%%%%%%%%%%%%%%%%%%%%%%%%%%%%%%%%%%%%%%%%%%
%% 
\subsection{Continual Lie algebras}
\label{continual}
In this subsection we recall the notion of a continual Lie algebra introduced in \cite{sv1}. 
It was then studied in \cite{sv2, sv3}.
 Suppose $E$ is an associative algebra 
over $\R$ or $\C$, and  
$K_0$, $K_{\pm 1}$, $K_{0,0}: E \times E \to E$,   
 are bilinear mappings.
 The local part of a continual Lie algebra can be defined as 
${\widehat{ \mathcal G}}=\mathcal G_{-1} \oplus \mathcal G_0 \oplus \mathcal G_{+1}$, 
 where $\mathcal G_i$, $i=0, \pm 1$.  
%% 
%%%%%%%%%%%%%%%%%%%%%%%%%%%%%%%%%%%%%%%%%%%%%%%%%%%%%%%%%%%%%%%%%% 
%%  
 The subspaces $\mathcal G_i$ consist of the elements
$\left\{ X_i(\phi), \phi \in E \right\}$, 
 $i=0, \pm 1$. 
The generators 
$X_i(\phi)$ are subject to commutation relations. 
%%
%%%%%%%%%%%%%%%%%%%%%%%%%%%%%%%%%%%%%%%%%%%%%%%%%%%%%%%%%%%%%%%%%%%%%
%% 
For instance, in the principle grading, 
they have the form 
%%
%%%%%%%%%%%%%%%%%%%%%%%%%%%%%%%%%%%%%%%%%%%%%%%%%%%%%%%%%%%%%%%%%%%%%%
%%
\begin{eqnarray*} 
&&\left[ X_0(\phi), X_0(\psi) \right]=X_0(K_{0,0}(\phi, \psi)),  
\nn
&&\left[ X_0(\phi), X_{\pm 1}(\psi) \right]=X_{\pm 1}(K_{\pm 1}(\phi, \psi)),
\nn
&&\qquad \left[ X_{+1}(\phi), X_{-1}(\psi) \right]=X_0(K_0(\phi, \psi)),
\end{eqnarray*}
%%
%%%%%%%%%%%%%%%%%%%%%%%%%%%%%%%%%%%%%%%%%%%%%%%%%%%%%%%%%%%%%%%%%%%%
%%
for all $\phi$, $\psi \in E$.
 It is also assumed that the Jacobi identity is satisfied.
 Then the conditions on mappings 
$\widehat K=\left(K_{0,0}, K_{0}, K_{\pm 1}\right)$ follow: 
%%
%%%%%%%%%%%%%%%%%%%%%%%%%%%%%%%%%%%%%%%%%%%%%%%%%%%%%%%%%%%%%%%%%%%
%%
\begin{eqnarray*} 
K_{\pm 1}(K_{0,0}( \phi, \psi) , \chi)&=&
K_{\pm 1}(\phi, K_{\pm 1}(\psi, \chi)) - K_{\pm 1}(\psi, K_{\pm 1}(\phi, \chi)),
\nn
K_{0,0}( \psi, K_0(\phi, \chi)) &=& 
K_0(K_{+1}(\psi, \phi), \chi)) + K_0(\phi, K_{-1}(\psi, \chi)),
\end{eqnarray*}
%%
%%%%%%%%%%%%%%%%%%%%%%%%%%%%%%%%%%%%%%%%%%%%%%%%%%%%%%%%%%%%%%%%%%%%
%%%
for all $\phi$, $\psi$, $\chi \in E$.
An infinite dimensional algebra 
$\mathcal G\left(E; \widehat K\right) = \mathcal G'\left(E; \widehat K\right)/J$,    
is called a continual contragredient Lie algebra,
 where $\mathcal G'(E; \widehat K)$ is a Lie algebra 
 freely generated by ${\widehat{ \mathcal G}}$, and $J$
is the largest homogeneous ideal with trivial intersection with $\mathcal G_0$ 
(consideration of the quotient is equivalent to imposing the Serre
relations in the case of finite-dimensional
simple complex Lie algebras) \cite{sv2, sv3}. 
%%
%%%%%%%%%%%%%%%%%%%%%%%%%%%%%%%%%%%%%%%%%%%%%%%%%%%%%%%%%%%%%%%%%%%%%%
%%
 $\mathcal G$ is endowed with a $\Z$-grading  
$\mathcal G= \bigoplus_{n \in \Z} \mathcal G_n$, 
where elements of subspaces $\mathcal G_n$ satisfy 
 the standard grading condition
$\left[ \mathcal G_n, \mathcal G_m \right] \subset \mathcal G_{n+m}$,  
where $X_n(\phi) \in \mathcal G_n$, 
and higher mapping $K_{n,m}(\phi, \psi)$ 
are present in commutation relations among $X_n(\phi)$ 
and $X_n(\psi)$.  
In general, the Jacobi identity for generators that belong to all grading spaces has the form  
\begin{eqnarray}
\left[X_i(\phi), \left[X_j(\psi), X_k(\theta) \right]  \right]
+ \left[X_j(\psi), \left[X_k(\theta), X_i(\phi) \right] \right] 
+ \left[X_k(\theta), \left[X_i(\phi), X_j(\psi) \right]  \right]  = 0,    
\nn
%%%%%%%%%%%%%%%%%%%%%%%%%%%%%%%%%%%%%%%%%%%%%%%%%%%%%%%%%%%%%%%%%%%%%%%%%%%%%  
%%
\label{jacjac}
 K_{i, j+k} \left(\phi, K_{j, k} \left(\psi, \theta\right) \right)   
+ K_{j, k+i } \left(\psi, K_{k, i} \left(\theta ,\phi \right) \right)   
+ K_{k, i+j} \left(\theta, K_{i, j} \left(\phi, \psi\right) \right) = 0.  
\end{eqnarray}
%% 
%%%%%%%%%%%%%%%%%%%%%%%%%%%%%%%%%%%%%%%%%%%%%%%%%%%%%%%%%%%%%%%%%%%%%%%%%%%%%%%%%%%%%%%%%%%%%
%%%%%%%%%%%%%%%%%%%%%%%%%%%%%%%%%%%%%%%%%%%%%%%%%%%%%%%%%%%%%%%%%%%%%%%%%%%%%%%%%%%%%%%%%%%%%%%%%%%
%%
\section{The algebra of differential relations from a chain complex}
\label{pseudo}
%%
%%%%%%%%%%%%%%%%%%%%%%%%%%%%%%%%%%%%%%%%%%%%%%%%%%%%%%%%%%%%%%%%%%%%%%%%%%%%%%%%%%%%%%%%%%%%%%%%
%% 
\subsection{The graded differential algebra}   
In this section we formulate the construction which allows us to generate an algebra of 
graded relations starting from a chain complex.  
%% 
%%%%%%%%%%%%%%%%%%%%%%%%%%%%%%%%%%%%%%%%%%%%%%%%%%%%%%%%%%%%%%%%%%%%%%%%%
%%
Consider a complex of spaces with elements depending 
on sets of parameters $\Theta_i$, $i \in \Z$, and given by  
\begin{equation}
\label{complex}
\ldots  
\stackrel{\delta(i-1)}{\longrightarrow} 
C(i, \Theta_i)   
\stackrel{\delta(i)}{\longrightarrow}C(i+1, \Theta_{i+1})    
 \stackrel{\delta(i+1)}{\longrightarrow} \ldots ,  
\end{equation}
%%
%%%%%%%%%%%%%%%%%%%%%%%%%%%%%%%%%%%%%%%%%%%%%%%%%%%%%%%%%%%%%%%%
%%
with a differential $\delta(i)$ satisfying the chain property  
$\delta(i+1) \circ \delta(i)=0$,   
for $i \in \Z$.   
%% 

%%%%%%%%%%%%%%%%%%%%%%%%%%%%%%%%%%%%%%%%%%%%%%%%%%%%%%%%%%%%%%%%%%%%%%%%%%%%%%%%%%%%%%%%%%%%%%
%%
We assume also that there exists a (not necessary associative)
 product among element of the spaces $C(i, \Theta_i)$,  
\begin{equation}
\label{initialproduct}
\cdot: C(i, \Theta_i) \cdot C(j, \Theta_j) \to C(f(i,j), \Theta_{g(i,j)}),   
\end{equation}
so that for any elements $\Phi_i \in C(i, \Theta_i)$ and  $\Phi_j \in C(j,\Theta_j)$,  
$\Phi_i \cdot \Phi_j=\Phi_{f(i,j)}$,  
where $f(i,j)$ and $g(i,j)$ are some functions of indices $i$ and $j$. 
Note that a function $f(i,j)$ defined the index 
for the product \eqref{initialproduct} resulting space. 
The function $g(i,j)$ defines the index for the resulting set of parameters.  
%%

%%%%%%%%%%%%%%%%%%%%%%%%%%%%%%%%%%%%%%%%%%%%%%%%%%%%%%%%%%%%%%%%%%%%%%%%%%%%%%%%%%%%%%%%
%%
Let us describe the functions $f(i,j)$ and $g(i,j)$ that we use in this paper explicitly. 
For the elements   
$\Phi_k \in   C(k, \Theta_k)$, and   
$\Phi_l \in C(l, \Theta_l)$, 
let $r$ be the number of common parameters among 
$k$- and $l$-sets of parameters 
 for $\Phi_k$ and $\Phi_l$. 
%%
%%%%%%%%%%%%%%%%%%%%%%%%%%%%%%%%%%%%%%%%%%%%%%%%%%%%%%%%%%%%%%%%%%%%%%%%%%%%%%%
%%
 We denote by $\deg \left(\Phi_k \right)=k$ the degree of an element  
 $\Phi_k \in C(k, \Theta_k)$. 
%%
%%%%%%%%%%%%%%%%%%%%%%%%%%%%%%%%%%%%%%%%%%%%%%%%%%%%%%%%%%%%%%%%%%%%%%%%%%%%%%%%%%%%%%
%%
In addition to the properties of a product $\cdot$ mentioned above, 
 we remove one set of $r$ common parameters from the product $\Phi_k \cdot \Phi_l$,  
 thus, the total degree of the product $\Phi_k \cdot \Phi_l$  
 is $\deg (\Phi_k \cdot \Phi_l)=k+l-r$.    
%%
%%%%%%%%%%%%%%%%%%%%%%%%%%%%%%%%%%%%%%%%%%%%%%%%%%%%%%%%%%%%%%%%%%%%%%%%%%%%%%%%%%%%%%
%% 
Let us assume that for all $i\in \Z$, the Leibniz rule formula  
for the operator $\delta(i)$ 
takes place, i.e.,   
%%
%%%%%%%%%%%%%%%%%%%%%%%%%%%%%%%%%%%%%%%%%%%%%%%%%%%%%%%%%%%%%%%%%%%%%%%%%%%%%%%%%%%%%%%%%%%
%% 
\begin{equation}
\label{leibniz0} 
  \delta(k+l) \left(\Phi_k \cdot \Phi_l \right)   
 = \left(\delta(k) \Phi_k \right) \cdot \Phi_l 
+ (-1)^{\deg \left(\Phi_k \cdot \Phi_l\right) } 
\Phi_k \cdot \delta (l) \Phi_l.   
\end{equation}    
%%
%%%%%%%%%%%%%%%%%%%%%%%%%%%%%%%%%%%%%%%%%%%%%%%%%%%%%%%%%%%%%%%%%%%%%%%%%%%%%%%%
%%
Note that in each part of the formula above, the total degree is $k+l-r+1$. 
%%
%%%%%%%%%%%%%%%%%%%%%%%%%%%%%%%%%%%%%%%%%%%%%%%%%%%%%%%%%%%%%%%%%%%%%%%%%%%%%%%%
%% 
We assume also that the chain complex above has some properties of an algebra 
with respect to the product \eqref{initialproduct}. 
%%
%%%%%%%%%%%%%%%%%%%%%%%%%%%%%%%%%%%%%%%%%%%%%%%%%%%%%%%%%%%%%%%%%%
%%
Let us assume, in particular, that the product \eqref{initialproduct} 
 satisfies the condition  
\begin{equation}
\label{prodcond}
\Phi_k  \cdot \Phi_l = - (-1)^{\deg \left(\Phi_k \cdot \Phi_l\right) }
\Phi_l \cdot \Phi_k. 
\end{equation}
for all $k$, $l \in \Z$,  $\Phi_k \in C(k, \Theta_k)$ and $\Phi_l \in C(l,\Theta_l)$. 
%%
%%%%%%%%%%%%%%%%%%%%%%%%%%%%%%%%%%%%%%%%%%%%%%%%%%%%%%%%%%%%%%%%%%%%%%%%%%%%%%%%%%%
%%
It is useful to introduce the following notation that we use in the text below:
 $\gamma_{n,m}=(-1)^{n m}$. 
%%
%%%%%%%%%%%%%%%%%%%%%%%%%%%%%%%%%%%%%%%%%%%%%%%%%%%%%%%%%%%%%%%%%%%%%%%%%%%%%%%%%%%%%%%%%%
%% 
Altogether, the complex \eqref{complex} form a graded differential algebra $\mathcal C$.   
In what follows, we skip the sets of parameters 
 $\Theta_i$ in the notations for the spaces of complexes, 
i.e., we set $C(i)= C(i, \Theta_i)$. 
Though one has to keep in mind that a product of elements of two 
spaces of complex \eqref{complex} depends on the resulting set of parameters 
$\Theta_{g(i, j)}$. 
%%
%%%%%%%%%%%%%%%%%%%%%%%%%%%%%%%%%%%%%%%%%%%%%%%%%%%%%%%%%%%%%%%%%%%%%%%%%%%%%%%%%%%%%%%%%%%%%
%% 
\subsection{The algebra of graded differential relations} 
%% 
%%%%%%%%%%%%%%%%%%%%%%%%%%%%%%%%%%%%%%%%%%%%%%%%%%%%%%%%%%%%%%%%%%%%%%%%%%%%%%%%%%%%%%%%%%%%%%
%%
As stated above, we assume that the product \eqref{initialproduct} 
is defined in such a way that common parameters of the 
sets $\Theta_i$ and $\Theta_j$ are present 
in the resulting set $\Theta_{g(i, j)}$ only once. 
%%
%%%%%%%%%%%%%%%%%%%%%%%%%%%%%%%%%%%%%%%%%%%%%%%%%%%%%%%%%%%%%%%%%
%%
For some $i$ and $j \in \Z$, 
let us introduce the orthogonality condition 
for a pair $C(i, \Theta_i)$ and  $C(j, \Theta_j)$,  
with respect to the product \eqref{initialproduct}. 
In particular, let us require that for a pair $C(i, \Theta_i)$, $C(j, \Theta_j)$,   
 there exist subspaces $C'(i) \subset C(i, \Theta_i)$ and $C'(j) \subset C(j,\Theta_j)$,   
 such that, for any $\Phi_i \in C'(i)$ and 
$\Phi_j \in C'(j)$,  
%%
%%%%%%%%%%%%%%%%%%%%%%%%%%%%%%%%%%%%%%%%%%%%%%%%%5
%%
\begin{equation}
\label{ortho}
\Phi_i \cdot \delta(j) \Phi_j=0, 
\end{equation}
%%
%%%%%%%%%%%%%%%%%%%%%%%%%%%%%%%%%%%%%%%%%%%%%%%%%
%%
 namely,  
$\Phi_i$ is supposed to be orthogonal to $\delta(j)\Phi_j$ 
with respect to \eqref{initialproduct}.    
By applying further differentials to \eqref{ortho} (and further consequences of such action),
 and using properties of
 a particular function $f(i,j)$ we obtain relations among elements of spaces 
$C(i, \Theta_i)$,
 $i\in \Z$.   
In particular, taking into account that both sides 
of such relations belong to the same space, 
we obtain limitations
(depending on the function $f$) on indices.  
In differential geometry, the orthogonality condition
 \eqref{ortho} provides the definition of 
integrability conditions 
for differential forms, and leads to the Frobenius theorem (see, e.g., \cite{Ghys}). 
%%

%%%%%%%%%%%%%%%%%%%%%%%%%%%%%%%%%%%%%%%%%%%%%%%%%%%%%%%%%%%%%%%%%%%
%%
Let us explain the notations we will use on few next pages. 
In \eqref{ogromno} we obtain (infinite) sequences 
of pairs of relations of the form 
\eqref{trojka}-\eqref{trojka1}. 
The (infinite) sequence of pair of 
relations has a tree graph structure with two sequences 
outgoing from one point. 
At each point we call one branch "left" and another branch 
as "right" marking the corresponding 
pair by $L$ or $R$. 
We denote such pairs by $\left(\alpha_i^{(K_i)}\right)$, 
where $\alpha_i^{(K_i)}$ is an element of
$C\left({n_i^{(K_i)}}\right)$ 
involved in 
 relations, and  $\left(K_i\right)$ is a sequence of $i$ entries each
  is either $L$ or $R$ for $i \ge 1$.    
%% 
 
%%%%%%%%%%%%%%%%%%%%%%%%%%%%%%%%%%%%%%%%%%%%%%%%%%%%%%%%%%%%%%%%%%%%%%%%%%%%%%%%%%%%%%%%%%%%%%%
%% 
Let $\chi \in C(n_0)$, $\Phi \in C(n)$, for any $n_0$, $n \in \Z$.      
%%
%%%%%%%%%%%%%%%%%%%%%%%%%%%%%%%%%%%%%%%%%%%%%%%%%%%%%%%%%%%%%%%%%%%%%%%%%%%%%
%%
Assume that for each $\Phi$ the annihilator
$\mathrm{Ann}_R(\Phi)=\left\{\phi\in C(i,\Theta_i)\mid\Phi\cdot\phi=0\right\}$  
is a principal right ideal, i.e., there exists an element $\alpha^{(R)}_1$
%%%
such that
$\mathrm{Ann}_R(\Phi) = \left\{\Phi \cdot C(j, \Theta_j)\right\}$. 
%% 
%%%%%%%%%%%%%%%%%%%%%%%%%%%%%%%%%%%%%%%%%%%%%%%%%%%%%%%%%%%%%%%%%%%%%%%
%% 
Then the orthogonality condition\eqref{ortho}  
applied to $\Phi$ and $\chi$ implies, i.e.,    
%% 
%%%%%%%%%%%%%%%%%%%%%%%%%%%%%%%%%%%%%%%%%%%%%%%%%%%%%%%%%%%%%%%%%%%%%%
%% 
\begin{equation}
\label{introrth0}
\Phi \cdot \delta(n_0) \chi=0,   
\end{equation}
%%
%%%%%%%%%%%%%%%%%%%%%%%%%%%%%%%%%%%%%%%%%%%%%%%%%%%%%%%%%%%%%%%%%
%% 
$\delta(n_0)\chi\in\mathrm{Ann}_R(\Phi)=\left\{\Phi\cdot C(k,\Theta_k)\right\}$, 
implies that there exists $\alpha^{(R)}_1\in C\left({n_1^{(R)}}\right)$ 
such that 
%% 
%%%%%%%%%%%%%%%%%%%%%%%%%%%%%%%%%%%%%%%%%%%%%%
%% 
\begin{equation}
\label{trojka}
 \delta(n_0) \chi= \Phi \cdot \alpha^{(R)}_1.  
\end{equation}
%%
%%%%%%%%%%%%%%%%%%%%%%%%%%%%%%%%%%%%%%%%%%%%%%%%%%%%%%%%%%%%%
%%
Let $r_1^{(R)}$ be the number of common parameters among 
$n$- and $n_1^{(R)}$-sets of parameters 
 for $\Phi$ and  $\alpha^{(R)}_1$. 
Since both sides of the last relation have to belong to the same space of the complex, 
the compatibility condition
\begin{equation}
\label{condi1}
n_0+1= n + n_1^{(R)} - r_1^{(R)},   
\end{equation}
occurs. 
%%
%%%%%%%%%%%%%%%%%%%%%%%%%%%%%%%%%%%%%%%%%%%%%%%%%%%%%%%%%%%%%%%
%%
Acting by $\delta(n_0+1)$ on \eqref{trojka} we obtain  
\begin{equation}
\label{trojka1}
 0 = \left(\delta(n) \Phi\right) \cdot \alpha^{(R)}_1 
+ (-1)^{(n_0+1)n_1^{(R)}}  
\Phi \cdot \delta\left(n^{(R)}_1 \right)\alpha^{(R)}_1.     
\end{equation}
%% 
%%%%%%%%%%%%%%%%%%%%%%%%%%%%%%%%%%%%%%%%%%%%%%%%%%%%%%%%%%%%%%%%%%%%%%%
%% 
Similarly, if the left annihilator
$\mathrm{Ann}_L(\delta(n_0)\chi)
=\left\{\psi\mid\psi\cdot\delta(n_0)\chi=0\right\}$  
is principal, then 
 \eqref{introrth0} implies that there exists 
$\alpha_1^{(L)} \in C\left(n_1^{(L)}\right)$,  
such that 
%%    
%%%%%%%%%%%%%%%%%%%%%%%%%%%%%%%%%%%%%%%%%%%%%%%%%%%%%%%5
%% 
\begin{eqnarray}
\label{vtorur0}
  \Phi =  \alpha_1^{(L)} \cdot \delta(n_0) \chi,  
\end{eqnarray}
with the condition
%%
%%%%%%%%%%%%%%%%%%%%%%%%%%%%%%%%%%%%%%%%%%%%%%%%%5
%% 
\begin{equation}
\label{condi2}
n = n^{(L)}_1 +  n_0 + 1 - r_1^{(L)}, 
\end{equation}
where $r_1^{(l)}$ is the number of common parameters among $n_1^{(L)}$ and $(n_0 +1)$ for 
   $\alpha_1^{(L)}$ and $\delta(n_0) \chi$.   
%%
%%%%%%%%%%%%%%%%%%%%%%%%%%%%%%%%%%%%%%%%%%%%%%%%%%%%%%%%%%%%%%%%%%%%%%%%%%%%%%%%%%%%%%%%%%%%%%%%
%% 
Consequently applying the corresponding $\delta$ operators to 
\eqref{introrth0}, \eqref{trojka} and \eqref{vtorur0}
we obtain the system of relations: 
%% 
%%%%%%%%%%%%%%%%%%%%%%%%%%%%%%%%%%%%%%%%%%%%%%%%%%%%%%%%%%%%%%%%%%%%%%%%
%% 
\begin{eqnarray}
\label{ogromno}
0&=&\Phi \cdot \delta(n_0) \chi,  \qquad (1) 
\nn
{\rm  \; (1) } \; \Rightarrow \; \Phi &=&  \alpha_1^{(L)} \cdot \delta(n_0) \chi, \qquad  (2)  
\nn
{  \delta. \rm (2)\; } \Rightarrow \; \delta(n)  
\Phi & = & \left( \delta\left(n^{(L)}_1\right)  \alpha_1^{(L)} \right) \cdot \delta(n_0) \chi,  \; \; (2')
 \;  \Rightarrow 0, 
\nn
{ \delta. \rm (1)} \; \Rightarrow \; 0&=& \left(\delta(n) \Phi \right) \cdot \delta(n_0) \chi, \qquad  (3)   
\nn
{\rm   (1)} \; \Rightarrow  \;  \delta(n_0) \chi&=& \Phi \cdot \alpha^{(R)}_1.
%%
%%%%%%%%%%%%%%%%%%%%%%%%%%%%%%%%%%%%%%%%%%%%%%%%%%%%%%%%%%%%%%%%%%%%%%%%%%%%%%%%%%% 
\nn
0&=&\delta(n^{(R)}_1) \Phi \cdot \alpha^{(R)}_1
+ (-1)^{n n^{(R)}_1}
 \Phi \cdot \delta\left(n^{(R)}_1\right) \alpha^{(R)}_1 \; \Rightarrow 0,   
\\ 
%% 
%%%%%%%%%%%%%%%%%%%%%%%%%%%%%%%%%%%%%%%%%%%%%%%%%%%%%%%%%%%%%%%%%%%%%%%%%%%%%%%%%
%% 
\nn
\label{ser1}
{\rm   (3) }\; \; \Rightarrow  \; \delta(n_0) \chi &=& \left(\delta(n) \Phi \right) \cdot \alpha^{(RR)}_2,    
\nn
  {\rm (4):} \; \;    0&=& \left( \delta(n) \Phi \right)\cdot \delta\left(n_2^{(RR)}\right) \alpha^{(RR)}_2   
 \Rightarrow \left(\alpha^{(RRR)}_3 \right)  
\Rightarrow ...
\nn
&  & \Downarrow  \qquad \qquad \qquad  \qquad \qquad \qquad  \qquad \qquad \qquad \qquad  \qquad (3')
\nn
&& \left( \alpha^{(RRL)}_3 \right) \Rightarrow ...  
\nn
&& \Downarrow
\nn
&& ... 
%%
%%%%%%%%%%%%%%%%%%%%%%%%%%%%%%%%%%%%%%%%%%%%%%%%%%%%%%%%%%%
%%
\\
\label{ser2}
{\rm   (3) }\; \Rightarrow\;  \delta(n) \Phi & = & \alpha^{(LL)}_2 \cdot \delta(n_0) \chi,    
\nn
{\rm (5):} \; \;  0 & = & \delta  \alpha^{(LL)}_2 \cdot \delta(n_0) \chi \Rightarrow \left( \alpha^{(LLR)}_3 \right)  
 \Rightarrow ...
\nn
&& \Downarrow \qquad \qquad \qquad  \qquad \qquad \qquad  \qquad \qquad \qquad \qquad  \qquad (3'')
\nn
&& \left( \alpha^{(LLL)}_3  \right) \Rightarrow ... 
\nn
&& \Downarrow  
\nn
&& ... 
\end{eqnarray} 
where we obtain (infinite) sequences  \eqref{ser1} and \eqref{ser2} of pairs of relations for 
$\alpha_i^{(K_i)} \in C\left({n_i^{(K_i)}}\right)$, $i \ge 1$. Recall that we denote such relations as  
$\left(\alpha_i^{(K_i)}  \right)$. 
%%
%%%%%%%%%%%%%%%%%%%%%%%%%%%%%%%%%%%%%%%%%%%%%%%%%%%%%%%%%%%%%%%%%%%%%%%%%%%%%%%%%%%%%%%%%%%%%%%%%%%%%%%%%%%%%%
%%
The corresponding indices  $n_i^{(K_i)}$  satisfy 
(in addition to \eqref{condi1} and \eqref{condi2}) relations   
for the sequence starting from $(4)$:
%%
%%%%%%%%%%%%%%%%%%%%%%%%%%%%%%%%%%%%%%%%%%%%%%%%%%%%%%%%%%%%%%%
%%
\begin{equation}
\label{rel2}
n_0=n + n_i^{(RRK_{i+2})} - r_i^{(RRK_{i+2})}, 
\end{equation}
for the sequence starting from $(5)$ for $i$, $j \ge 2$ 
\begin{equation}
\label{rel3}
n=n_0 + n_j^{(LLK_{j+2})} - r_j^{(LLK_{j+2})}. 
\end{equation}
%%
%%%%%%%%%%%%%%%%%%%%%%%%%%%%%%%%%%%%%%%%%%%%%%%%%%%%%%%%%%%%%%%%%%%%%%%%%%%%%%%%%%%%%%%%%%%%%%%%
One can easily see that not all elements in \eqref{ogromno} are independent. 
For instance, from \eqref{condi1} and \eqref{condi2} we obtain 
\[
\left(n^{(L)}_1 - r^{(L)}_1\right)  = - \left( n_1^{(R)}- r^{(R)}_1 \right) = n - n_0 -1.
\]
%% 
%%%%%%%%%%%%%%%%%%%%%%%%%%%%%%%%%%%%%%%%%%%%%%%%%%%%%%%%%%%%%%%%%%%%%%%%%%%    
%%
From \eqref{trojka}-\eqref{vtorur0}, and from $(3')$-$(3'')$  of \eqref{ogromno} 
we infer that $\alpha_1^{(L)}$, $\alpha_1^{(R)}$, and $\alpha_i^{(LLK_i)}$, $\alpha_i^{(RRK_i)}$ 
are related by a conjugation with respect to the product \eqref{initialproduct}: 
%%
%%%%%%%%%%%%%%%%%%%%%%%%%%%%%%%%%%%%%%%%%%%%%%%%%%%%%%%%%%%%%%%%%%%%%%%%%%%%
%% 
\begin{eqnarray*}
&&\Phi= \alpha_1^{(L)} \cdot \left( \Phi \cdot  \alpha_1^{(R)}\right),    
\nn
&&\delta(n_0) \chi= \left( \alpha_i^{(LLK_i)}
 \cdot \delta(n_0)\; \chi \right) \cdot  \alpha_i^{(RRK_i)}.
\end{eqnarray*}        
%%
%%%%%%%%%%%%%%%%%%%%%%%%%%%%%%%%%%%%%%%%%%%%%%%%%%%%%%%%%%%%%%%%%%%%%%%%
%% 
Similar relations apply among other elements $\alpha_i^{(K_i)}$.  
%%

%%%%%%%%%%%%%%%%%%%%%%%%%%%%%%%%%%%%%%%%%%%%%%%%%%%%%%%%%%%%%
%% 
%%
We supply here a recursive representation summarizing the systems 
\eqref{ogromno}-\eqref{ser2} in the following form. 
%% 
%%%%%%%%%%%%%%%%%%%%%%%%%%%%%%%%%%%%%%%%%%%%%%%%%%%%%%%%%%%%%%%%%%%%%%%%%
%%
In order to establish a clear way to rewrite such a system
 we make the recursion explicit over the word $K_i$. 
%% 
%%%%%%%%%%%%%%%%%%%%%%%%%%%%%%%%%%%%%%%%%%%%%%%%%%%%%%%%%%%%%%%%%%%%%
%% 
The indexing set consists of the following. 
Let $\mathcal{K}= \{L,R\}^*$ be the set of all finite words 
over $\{L,R\}$, and let $\epsilon$ denote the empty word.
%%
%%%%%%%%%%%%%%%%%%%%%%%%%%%%%%%%%%%%%%%%%%%%%%%%%%%%%%%%%%%%%%%%%%%%%
%% 
The generators are the following. 
For each $i \ge 1$, define
%% 
%%%%%%%%%%%%%%%%%%%%%%%%%%%%%%%%%%%%%%%%%%%%%%%%%%%%%%%%%%%%%%%
\[
\alpha_i : \mathcal{K} \to \bigsqcup_{n \in \mathbb{Z}} C(n), 
\qquad
n_i : \mathcal{K} \to \mathbb{Z}, 
\qquad
r_i : \mathcal{K} \to \mathbb{Z}_{\ge 0}.
\]
%%
%%%%%%%%%%%%%%%%%%%%%%%%%%%%%%%%%%%%%%%%%%%%%%%%%%%%%%%%%%%%%
%%
Now we proceed with levels of the system \eqref{ogromno}-\eqref{ser2}. 
The first level ($i=1$) is determined by \eqref{trojka} 
 and \eqref{vtorur0}:  
%%
%%%%%%%%%%%%%%%%%%%%%%%%%%%%%%%%%%%%%%%%%%%%%%%%%%%%%% 
\[
\delta(n_0)\chi = \Phi \cdot \alpha_1^{(R)}, 
\qquad
\Phi = \alpha_1^{(L)} \cdot \delta(n_0)\chi.
\]
%%%%%%%%%%%%%%%%%%%%%%%%%%%%%%%%%%%%%%%%%%%%%%%%%%%%%%
%% 
The indices satisfy the compatibility conditions:
\[
n_0 + 1 = n + n_1^{(R)} - r_1^{(R)}, 
\qquad
n = n_1^{(L)} + n_0 + 1 - r_1^{(L)}.
\]
%% 
%%%%%%%%%%%%%%%%%%%%%%%%%%%%%%%%%%%%%%%%%%%%%%%%%%%%%%%%
%%
The recursion: let $K \in \mathcal{K}$. Then we have 
 the R-extension: 
if $\alpha_i^{(K)}$ appears in a relation of the form
%%
%%%%%%%%%%%%%%%%%%%%%%%%%%%%%%%%%%%%%%%%%%%%%%%%%%%%%%
%%  
\[
\delta(n_0)\chi = (\cdots) \cdot \alpha_i^{(K)},
\]
%%
%%%%%%%%%%%%%%%%%%%%%%%%%%%%%%%%%%%%%%%%%%%%%%%%%%%%%%%%
%%  
then applying $\delta(\cdots)$ and using 
the graded Leibniz rule \eqref{leibniz0} 
 gives 
%%
%%%%%%%%%%%%%%%%%%%%%%%%%%%%%%%%%%%%%%%%%%%%%%%%%%%%%% 
%% 
\[
0 = (\cdots)' \cdot \alpha_i^{(K)} + (-1)^{\deg(\cdots) \deg(\alpha_i^{(K)})} 
(\cdots) \cdot \delta\bigl(n_i^{(K)}\bigr)\alpha_i^{(K)}.
\]
%%
%%%%%%%%%%%%%%%%%%%%%%%%%%%%%%%%%%%%%%%%%%%%%%%%%%%%%%%%%%
%% 
Hence there exists $\alpha_{i+1}^{(KR)}$ such that
\[
\delta(n_0)\chi = (\cdots)' \cdot \alpha_{i+1}^{(KR)}.
\]
%% 
%%%%%%%%%%%%%%%%%%%%%%%%%%%%%%%%%%%%%%%%%%%%%%%%%%%%%%%%
%% 
Thus 
$\alpha_{i+1}^{(KR)}$ is defined by 
$\delta(n_0)\chi = (\delta(\cdots)) \cdot \alpha_{i+1}^{(KR)}$. 
%% 
%%%%%%%%%%%%%%%%%%%%%%%%%%%%%%%%%%%%%%%%%%%%%%%%%%%%%%%%%%%%%%
%% 
The indices satisfy \eqref{rel2} 
\[
n_0 = n + n_{i+1}^{(KR)} - r_{i+1}^{(KR)}.
\]
%%
%%%%%%%%%%%%%%%%%%%%%%%%%%%%%%%%%%%%%%%%%%%%%%%%%%%%
%%
We have also the L-extension. 
If $\alpha_i^{(K)}$ appears in a relation of the form
%%
%%%%%%%%%%%%%%%%%%%%%%%%%%%%%%%%%%%%%%%%%%%%%%%%%% 
\[
\delta(n)\Phi = \alpha_i^{(K)} \cdot \delta(n_0)\chi,
\]
then applying $\delta(\cdots)$ gives
%%
%%%%%%%%%%%%%%%%%%%%%%%%%%%%%%%%%%%%%%%%%%%%%%%%%%%%% 
%% 
\[
0 = \delta\bigl(n_i^{(K)}\bigr)\alpha_i^{(K)} \cdot \delta(n_0)\chi.
\]
%% 
%%%%%%%%%%%%%%%%%%%%%%%%%%%%%%%%%%%%%%%%%%%%%%%%%%%%
%% 
Hence there exists $\alpha_{i+1}^{(KL)}$ such that
\[
\delta(n)\Phi = \alpha_{i+1}^{(KL)} \cdot \delta(n_0)\chi.
\]
%% 
%%%%%%%%%%%%%%%%%%%%%%%%%%%%%%%%%%%%%%%%%%%%%%%
Thus, 
$\alpha_{i+1}^{(KL)}$ is defined by 
$\delta(n)\Phi = \alpha_{i+1}^{(KL)} \cdot \delta(n_0)\chi$. 
%%% 
%%%%%%%%%%%%%%%%%%%%%%%%%%%%%%%%%%%%%%%%%%%% 
%% 
The indices satisfy \eqref{rel3}  
\[
n = n_0 + n_{i+1}^{(KL)} - r_{i+1}^{(KL)}.
\]
%%
%%%%%%%%%%%%%%%%%%%%%%%%%%%%%%%%%%%%%%%%%%%%%%%%%%%%%%%%%%
%% 
The explicit recursive formulation: 
for all $K \in \mathcal{K}$, 
$\alpha_{i+1}^{(KR)}$ is the unique element such that 
\[
\delta(n_0)\chi = (\delta^{(\cdots)}\Phi)\cdot \alpha_{i+1}^{(KR)},
\] 
$\alpha_{i+1}^{(KL)}$ is the unique element such that  
\[
\delta(n)\Phi = \alpha_{i+1}^{(KL)} \cdot \delta(n_0)\chi.
\]
%%
%%%%%%%%%%%%%%%%%%%%%%%%%%%%%%%%%%%%%%%%%%%%%%%%%%%%%%%%%%%%%
%% 
The compact operator form of the system is the following: 
define maps  $T_R(\alpha)=\alpha'$ such that 
\[
\delta(n_0)\chi=(\delta\Phi)\cdot \alpha',
\]
and $T_L(\alpha)=\alpha'$ such that 
\[
\delta(n)\Phi = \alpha' \cdot \delta(n_0)\chi.
\]
%%
%%%%%%%%%%%%%%%%%%%%%%%%%%%%%%%%%%%%%%%%%%%%%%%%%%%
%% 
Then for $K = k_1 \cdots k_m$,
\[
\alpha_i^{(K)} = T_{k_m} \circ \cdots \circ T_{k_1}(\alpha_1^{(\cdot)}).
\]
%%
%%%%%%%%%%%%%%%%%%%%%%%%%%%%%%%%%%%%%%%%%%%%%%%%%%%%%%%%%%%%%%%%%%%%%%%%%%%%%%%%%%%
%%
All that can be summarized in the table: 
%%
%%%%%%%%%%%%%%%%%%%%%%%%%%%%%%%%%%%%%%%%%%%%%%%%%%%%%%%%%%%%%%%%%%%%%%%%%%%%%%%
%%
\[
\begin{array}{|c|c|c|c|}
\hline
\textbf{Level} & \textbf{Branch} & \textbf{Relation} & \textbf{Recursive form} \\
\hline

0 & - 
& 0 = \Phi \cdot \delta(n_0)\chi 
& \mathcal{R}_0= \Phi \perp \delta(n_0)\chi
\\
\hline

1 & R 
& \delta(n_0)\chi = \Phi \cdot \alpha^{(R)}_1 
& \alpha^{(R)}_1 \in C(n^{(R)}_1)
\\

1 & L 
& \Phi = \alpha^{(L)}_1 \cdot \delta(n_0)\chi 
& \alpha^{(L)}_1 \in C(n^{(L)}_1)
\\
\hline

2 & RR 
& \delta(n_0)\chi = (\delta(n)\Phi)\cdot \alpha^{(RR)}_2 
& \alpha^{(RR)}_2 = \mathcal{F}_R(\alpha^{(R)}_1)
\\

2 & LL 
& \delta(n)\Phi = \alpha^{(LL)}_2 \cdot \delta(n_0)\chi 
& \alpha^{(LL)}_2 = \mathcal{F}_L(\alpha^{(L)}_1)
\\
\hline

i \ge 1 & R K_i 
& \delta(\cdot) = (\delta(\cdot)) \cdot \alpha^{(R K_i)}_{i+1} 
& \alpha^{(R K_i)}_{i+1} = \mathcal{F}_R\big(\alpha^{(K_i)}_i\big)
\\

i \ge 1 & L K_i 
& \delta(\cdot) = \alpha^{(L K_i)}_{i+1} \cdot \delta(\cdot) 
& \alpha^{(L K_i)}_{i+1} = \mathcal{F}_L\big(\alpha^{(K_i)}_i\big)
\\
\hline

\end{array}
\]
%%%%%%%%%%%%%%%%%%%%%%%%%%%%%%%%%%%%%%%%%%%%%%%%%%%%%%%%%%%%%%%%55
%% 
\[
\alpha^{(K_{i+1})}_{i+1} =
\begin{cases}
\mathcal{F}_R\big(\alpha^{(K_i)}_i\big), & K_{i+1} = (K_i, R), 
\\[6pt]
\mathcal{F}_L\big(\alpha^{(K_i)}_i\big), & K_{i+1} = (K_i, L),
\end{cases}
\]
%%%%%%%%%%%%%%%%%%%%%%%%%%%%%%%%%%%%%%%%%%%%%%%%%%%%%%%%%%
%%
with induced relations:
 the right branch: $\delta(\cdot) = (\delta(\cdot)) \cdot \alpha^{(K_i)}_i$,  
 the left branch:  $\delta(\cdot) = \alpha^{(K_i)}_i \cdot \delta(\cdot)$.  
%% 
%%%%%%%%%%%%%%%%%%%%%%%%%%%%%%%%%%%%%%%%%%%%%%%%%%%%%%%%%%%%%%%%%%%%%%%%%%%%%%%%%%%%
%% 
The system forms a binary tree of relations. 
Each node corresponds to an element
$
\alpha_
i
^{(K_i)}$. 
Paths $
K_
i=(L, R)$ encode iterated left/right factorizations. 
(2.13) \eqref{ser1} corresponds to the 
right subtree
, $RRK_i$ 
(2.14) \eqref{ser2} to the left subtree 
$
LLK_i$. 
%% 
%%

%%%%%%%%%%%%%%%%%%%%%%%%%%%%%%%%%%%%%%%%%%%%%%%%%%%%%%%%%%%%%%%%%%%%%%%
%%
The sequence of relations \eqref{ogromno} cancels when one of relations \eqref{condi1}-\eqref{condi2} or 
\eqref{rel2}-\eqref{rel3} for a sequence of pairs of equations is not fulfilled.  
The natural grading is given by the condition that both sides of 
 relations in \eqref{ogromno} belong to the same chain complex space. 
%%
%%%%%%%%%%%%%%%%%%%%%%%%%%%%%%%%%%%%%%%%%%%%%%%%%%%%%%%%%%%%%%%%%%%%%%%%%%%%%%%%%%%%%%%%%%%%%%%
%%
In this paper we consider the simplest form \eqref{complex} of a complex.
In general, for more complicated actions of $\delta$, such that 
 $\delta(i): C(i) \to   C(i+k(i))$,    
where $k(i)$ depends on $i\in \Z$ 
(see, e.g., \cite{Huang} for non-trivial actions of certain $\delta$ operators among chain 
complex spaces 
for vertex algebras. 
Consideration of such complexes will be given by the author in another article).  
The corresponding  relations as well 
as compatibility relations  could be different 
from \eqref{ogromno}
%%

%%%%%%%%%%%%%%%%%%%%%%%%%%%%%%%%%%%%%%%%%%%%%%%%%%%%%%%%%%%%%%%%%%%%%%%%%%
As an upshot, 
the orthogonality condition \eqref{ortho}  for all choices of $n_0$, $n \in \Z$,   
and the conditions \eqref{rel2}-\eqref{rel3},  
 applied to 
the chain complex \eqref{complex}  
bring about the structure 
of a graded algebra with  relations \eqref{ogromno} with   
 respect to the multiplication \eqref{prodcond}.   
%%
%%%%%%%%%%%%%%%%%%%%%%%%%%%%%%%%%%%%%%%%%%%%%%%%%%%%%%%%%%%%%%%%%%%%%%%%%%%%%%%%%%%%%%%%%%%%%%
%%
%%
As we can see, the system of relations \eqref{ogromno} has a tree structure.  
"Left" and "right" directions have a mixture of dependent elements.  
Let us denote $n(i)= n_i^{(K_i)}$. 
For $n_0$, $n\in \Z$, let $I(n_0, n)$ be the set of sequences of indices $(n(i))$, $i\ge 0$,  
marking all paths $(K_i)$ 
in the tree structure of \eqref{ogromno},
 describing "left" or "right" choice at each point. 
Let 
$I= \bigcup\limits_{n_0, n \in \Z} \;  I(n_0, n)$,  
be the space of all sequences over the tree graph for a complex \eqref{complex}. 
Denote by $I_0 \subset I$ the subset of such paths that include pairs of independent elements only.  
Then we are able to single out generators and commutation relations of   
a Lie algebra with a continual space of roots,  
 (see subsection \eqref{continual} and \cite{sv1}).     
%%  
%%%%%%%%%%%%%%%%%%%%%%%%%%%%%%%%%%%%%%%%%%%%%%%%%%%%%%%%%%%%%%%%%%%%%%%%%%%%%%%%%%%%%%%%%%%%%%%%
%% 
\subsection{Construction of a continual Lie algebra from a chain complex}
Let us further assume that the spaces $C(i)$, $i \in \Z$ 
admit an ordinary (not necessary commutative)   
product $\Phi \Psi$ among elements for $\Phi_k \in C(k)$,  
and $\Phi_l \in C(l)$, for all $k$, $l \in \Z$. 
%%
%%%%%%%%%%%%%%%%%%%%%%%%%%%%%%%%%%%%%%%%%%%%%%%%%%%%%%%%%%%%%%%
%%
Then, as a product \eqref{initialproduct}, 
satisfying conditions \eqref{leibniz0} and \eqref{prodcond}, one can take     
%%
%%%%%%%%%%%%%%%%%%%%%%%%%%%%%%%%%%%%%%%%%%%%%%%%%%%%%%%%%%%%%%%%
%%
 \begin{eqnarray}
\label{product0}
 \Phi_k \cdot \Phi_l = \left[\Phi_k, \Phi_l\right]_{gr} 
=\Phi_k\Phi_l-(-1)^{\deg \left(\Phi_k \cdot \Phi_l\right)}\Phi_l \Phi_k,      
\end{eqnarray}
where brackets mean the ordinary graded commutator.   
%%
%%%%%%%%%%%%%%%%%%%%%%%%%%%%%%%%%%%%%%%%%%%%%%%%%%%%%%%%%%%%%%%%%
%% 
It is known that the introduction of the commutator 
with respect to the original multiplication of an algebra 
transfers it into a Lie algebra when the Jacobi conditions are satisfied. 
In our case we show that the differential algebra $\mathcal C$ 
defined above being supplied with 
the orthogonality conditions deliver the structure of a continual Lie algebra. 
%% 
%%%%%%%%%%%%%%%%%%%%%%%%%%%%%%%%%%%%%%%%%%%%%%%%%%%%%%%%%%%%%%%%%%%%%%%%%%
%%
The setup of this subsection combined with the previous subsection   
provides us with a proof of the main result of this paper, Theorem \ref{mainpro}. 
%%

%%%%%%%%%%%%%%%%%%%%%%%%%%%%%%%%%%%%%%%%%%%%%%%%%%%%%%%%%%%%%%%%%%%%%%%%%%%%%%%%%%%%%%%%%%%%%%
%%
Though the structure of the system \eqref{ogromno} may seem to be not very complicated, 
actual properties of the corresponding continual Lie algebra depend on properties of 
the spaces of a specific the bicomplex \eqref{complex} and the nature of parameters $\Theta_i$. 
%%
%%%%%%%%%%%%%%%%%%%%%%%%%%%%%%%%%%%%%%%%%%%%%%%%%%%%%%%%%%%%%%%%%%%%%%%%%%%%%%%%%%%%
%% 
For a fixed choice of a path of independent functions/differential relations in the system 
\eqref{ogromno}, 
there exists a variety of choices on how to identify generators 
of a continual Lie algebra with 
generators of the differential algebra.  
Therefore, the actual form of commutation relations for the corresponding 
continual Lie algebra varies accordingly. 
One can also chose various ways how to define 
a grading for each specific $\mathcal G(\Theta_{n(i)})$, 
for generators of a continual Lie algebra
 (see subsection \ref{continual}) resulting from 
\eqref{ogromno}. 
%%
%%%%%%%%%%%%%%%%%%%%%%%%%%%%%%%%%%%%%%%%%%%%%%%%%%%%%%%%%%%%%%
%%
The structure of the product \eqref{initialproduct},  
 the condition \eqref{leibniz0}, and the action of 
the differentials $\delta$ provide  
  the Jacobi identity for generators on the continual Lie algebra 
$\mathcal G(\Theta_{n(i)})$, and, simultaneously, apply conditions of the form \eqref{jacjac}
to elements of the parameter spaces $\Theta_{(n(i))}$, $i \ge 0$.      
In the next Section we specify the above construction and Theorem \eqref{mainpro}  
 in the case of  double cochain complex \eqref{Cpq}-\eqref{deltas} associated with  
%%%
differential forms \cite{CM}. We derive explicitly the generators and commutation relations
 for the corresponding continual Lie algebras.  
%%
%%%%%%%%%%%%%%%%%%%%%%%%%%%%%%%%%%%%%%%%%%%%%%%%%%%%%%%%%%%%%%%%%%%%%%%%%%%%%%%%%%%
%%%%%%%%%%%%%%%%%%%%%%%%%%%%%%%%%%%%%%%%%%%%%%%%%%%%%%%%%%%%%%%%%%%%%%%%%%%%%%%%%%%
%% 
\section{An example: double complex associated with foliations} 
\label{example}
%%
%%%%%%%%%%%%%%%%%%%%%%%%%%%%%%%%%%%%%%%%%%%%%%%%%%%%%%%%%%%%%%%%%%%%%%%%%%%%%%%%%%%%%%%%%%
%% 
\subsection{\v{C}ech-de Rham complex for foliations}
\label{cm}
%%
%%%%%%%%%%%%%%%%%%%%%%%%%%%%%%%%%%%%%%%%%%%%%%%%%%%%%%%%%%%%%%%%%%%%%%%%%%%%%%%%%%%%%%%%%%%%%
%%
In this subsection we recall \cite{CM} the notion of the basis of transversal sections,  
and \v{C}ech-de Rham complex for a foliation of a smooth manifold.  
 Let $M$ be a smooth manifold of dimension $n$, 
equipped with a foliation $\F$ of co-dimension $l$ \cite{CM}. 
%%
%%%%%%%%%%%%%%%%%%%%%%%%%%%%%%%%%%%%%%%%%%%%%%%%%%%%%%%%%%%%%%%
%%
A transversal section of $\F$ is
 an embedded $l$-dimensional submanifold $U\subset M$ which 
is everywhere transverse to $\F$ leaves. 
If $\alpha$ is a path between two points 
$x$ and $y$ on the same leaf, and if $U$ and $V$ are transversal sections through $x$ and $y$,
 then $\alpha$ defines a transport
along the leaves from a neighborhood of $x$ in $U$ to a neighborhood of $y$ in $V$. 
%%
%%%%%%%%%%%%%%%%%%%%%%%%%%%%%%%%%%%%%%%%%%%%%%%%%%%%%%%%%%%%%%%%%%%%%%%%%%%%%%%%%%%%
%%
Therefore, we define a germ of a diffeomorphism $hol(\alpha): (U, x)\rmap (V, y)$,
called the holonomy of the path $\alpha$. 
%%%
%%%%%%%%%%%%%%%%%%%%%%%%%%%%%%%%%%%%%%%%%%%%%%%%%%%%%
%%
If such a transport
is defined in all of $U$ and embeds $U$ into $V$, this embedding 
$h: U\hookrightarrow V$ is sometimes also denoted by $hol(\alpha): U\hookrightarrow V$.  
Embeddings of this form will be called holonomy embeddings.  
%%

%%%%%%%%%%%%%%%%%%%%%%%%%%%%%%%%%%%%%%%%%%%%%%%%%%%%%%%%%
%%
Transversal sections $U$ through a point $x$ are 
  neighborhoods of the leaf through $x$ in the leaf space. 
%%
%%%%%%%%%%%%%%%%%%%%%%%%%%%%%%%%%%%%%%%%%%%%%%%%%%%%%%%%%%%%%%
 One defines  a  transversal basis for $(M, \F)$ as a 
family $\U$ of transversal sections $U\subset M$ 
with the property that,
 if $V$ is any transversal section through a given point $y\in M$,
 there exists a 
holonomy embedding $h: U\hookrightarrow V$ with $U\in \U$ and $y\in h(U)$.
%%
%%%%%%%%%%%%%%%%%%%%%%%%%%%%%%%%%%%%%%%%%%%%%%%%%%%%%%%%%%%%%%%%%%%%%%%%%%%%
%%
  A transversal section is a $l$-disk given by a chart for the foliation. 
One then constructs   
a transversal basis $\U$ out of a basis $\tilde{\U}$ of $M$ by domains of foliation charts
 $\phi_{U}: \tilde{U}\tilde{\rmap} \mathbb{R}^{n-l}\times U$, $\tilde{U}\in \tilde{\U}$,
 with $U=\mathbb{R}^l$.
%%
%%%%%%%%%%%%%%%%%%%%%%%%%%%%%%%%%%%%%%%%%%%%%%%%%%%%%%%%%%%%%%%%%%%%%%%%%%%%%%%%%%%%%%%%%%%%%%

%%
Let us recall the construction of the  \v{C}ech-de~Rham 
cohomology in \cite{CM}.
Let  $\U$ be a family of transversal sections of $\F$. 
%%
%%%%%%%%%%%%%%%%%%%%%%%%%%%%%%%%%%%%%%%%%%%%%%%%%%%%%%%%%%%%%%%%%
%%
 Consider the double complex 
\begin{equation}
\label{Cpq}
C^{p,q}(\F)=\prod_{U_0\stackrel{h_1}{\longrightarrow}\cdots
\stackrel{h_p}{\longrightarrow} U_p}\Omega^q(U_0),
\end{equation}
%%
%%%%%%%%%%%%%%%%%%%%%%%%%%%%%%%%%%%%%%%%%%%%%%%%%%%%%%%%%%%%%
%%
 where the product ranges over all $p$-tuples 
of holonomy embeddings $h_i$, $0\le i \le p$, between  
transversal sections from a fixed transversal basis $\U$, and $\Omega^q$ 
is the space of differential forms of order $q$. 
%%
%%%%%%%%%%%%%%%%%%%%%%%%%%%%%%%%%%%%%%%%%%%%%%%%%%%%%%%%%%%%%
%%
The vertical differential is defined as 
$(-1)^p d:C^{p,q}(\F)\to C^{p,q+1}(\F)$,  
 where $d$ is the usual de Rham differential. 
The horizontal differential 
$\mathcal D:C^{p,q}(\F) \to C^{p+1,q}(\F)$,  
 is given by
$\mathcal D= \sum\limits_{i=0}^{k+1}(-1)^i\mathcal D_i$,  
where
\begin{equation}
\label{deltas} 
\mathcal D_i \omega(h_1, \ldots , h_{k+1})= \left\{
 \begin{array}{lll}
 h_1^*\omega(h_2, \ldots , h_{k+1}), \; \mbox{if $i=0,$}\\ 
\omega(h_1, \ldots, h_{i+1}h_i, \ldots, h_{k+1}), \; \mbox{if $0<i< k+1,$}\\
 \omega(h_1, \ldots, h_k),  \; \mbox{if $i= k+1.$}
  \end{array}
 \right.
\end{equation}
%%
%%%%%%%%%%%%%%%%%%%%%%%%%%%%%%%%%%%%%%%%%%%%%%%%%%%%%%%%%%%%%%%%%%%%%%%%
%%
 The product is defined by 
%%
%%%%%%%%%%%%%%%%%%%%%%%%%%%%%%%%%%%%%%%%%%%%%%%%%%%%%%%%%%%%%%%%
%%
\begin{equation}
\label{bigradif}
(\omega \eta)\left(h_1, \ldots , h_{n+n'}\right)= (-1)^{nn'} 
\omega(h_1, \ldots , 
h_n) \; \left( h_1^* \ldots h_n^* \right).\eta\left(h_{n+1}, \ldots, h_{n+n'}\right),     
\end{equation}
for $\omega\in C^{n, m}$ and $\eta\in C^{n',m'}$,
 and $h^*_i$ being the dual to $h_i$.  
Thus $(\omega\;\eta)(h_1, \ldots , h_{n+n'}) \in C^{n+ n',m+ m'}$, 
and the product  \eqref{bigradif}  
 delivers the structure of a bigraded differential algebra. 
The cohomology of this complex is called the \v{C}ech-de~Rham
cohomology $\vH^*_\U(M/\F)$ of the leaf space $M/\F$ with respect 
to the transversal basis $\U$.  
%%
%%%%%%%%%%%%%%%%%%%%%%%%%%%%%%%%%%%%%%%%%%%%%%%%%%%%%%%%%%%%%%%%%%%%%%%%%%%%%%%%%%
%% 
\subsection{The double complex} 
In the previous subsection we recalled 
 the formulation of the \v{C}ech-de~Rham cohomology given in \cite{CM} 
for foliations on smooth manifolds (see subsection \ref{cm}).    
In this case the spaces in \eqref{complex} 
are differential forms $C^{n,m}(\F)$ defined on a foliation,     
and the coboundary operators is given by 
$\delta^{p, q} = (-1)^p d + \delta^q$, 
 where $d$ and $\delta$ operators are
defined in subsection \ref{cm}.   
The ordinary product for differential forms
 $\omega^{n,m}(h_1, \ldots, h_m) \in C^{n,m}\left(\F \right)$ is    
given by \eqref{bigradif}.  
 The product \eqref{initialproduct} required 
for the formulation of Section \ref{pseudo} is provided by 
the commutator \eqref{product0}.  
As it was mentioned in \cite{CM}, the bicomplex \eqref{Cpq}-\eqref{deltas} has the 
structure of a 
 bigraded differential algebra
with respect to the ordinary product \eqref{bigradif}. 
%%
%%%%%%%%%%%%%%%%%%%%%%%%%%%%%%%%%%%%%%%%%%%%%%%%%%%%%%%%%%%%%%%
%%
According to explanations of Section \ref{pseudo},  
for the bicomplex \eqref{Cpq}-\eqref{deltas} 
generates a graded differential algebra with relations  
\eqref{ogromno} with respect to the product \eqref{product0}.  
%%

%%%%%%%%%%%%%%%%%%%%%%%%%%%%%%%%%%%%%%%%%%%%%%%%%%%%%%%%%%%%%%%%%%%%%%%%%%%%%%%%%%%%%%%%%%%%%%
%%
 Let us recall that we assume 
 non-negative 
indices for all bicomplex spaces $C^{n_i,m_i}$.   
Let  $\chi \in C^{n_0, m_0} $, $\Phi \in C^{n, m}$.  
The orthogonality condition \eqref{ortho} with respect ot the product \eqref{product0}
leads to systems of the form \eqref{ogromno} when
 applied to the double complex \eqref{complex}. 
%% 
%%%%%%%%%%%%%%%%%%%%%%%%%%%%%%%%%%%%%%%%%%%%%%%%%%%%%%%%%%%%%%%%%%%%%%%%%%%%%
%% 
In particular, for $\chi \in C^{n_0,m_0} $ and $\Phi \in C^{n_0,m_0}$,      
we obtain the system \eqref{ogromno} of relations for elements 
$\chi$, $\delta^{n_0,m_0}\chi$, $\Phi$, $\delta^{n,m} \Phi$,  
$\alpha_i^{(K_i)} \in C^{n_i^{(K_i)},m_i^{(K_i)}} $,   
$\delta^{n(i),m(i)} 
 \alpha_i^{(K_i)}$, $i \ge 1$, for $n\ge 0$, $m \ge 0$.  
Let $r^{(K_i)}_i$ and $t^{(K_i)}_i$ be numbers of common degrees    
and transversal sections of $C^{n_i,m_i}$
 for the forms $\chi \in C^{n_0,m_0} $ and $\Phi \in C^{n_0,m_0}$.  
Then 
the compatibility conditions \eqref{condi1}-\eqref{condi2} 
and \eqref{rel2}-\eqref{rel3} for 
indices  
$n$, $m$, $n_0$, $m_0$, $n_i^{(K_i)}$, $m_i^{(K_i)}$, satisfy the relations in vector form:  
\begin{eqnarray}
\label{condition1}
(n_0 +1, m_0 +1) &=&(n, m)+ \left(n_1^{(R)},m_1^{(R)}\right) 
 - \left(r_1^{(R)}, t_1^{(R)}\right),     
\end{eqnarray}
for $(2)$ in \eqref{ogromno}:
\begin{eqnarray}
\label{dopcond1}
(n,m) &=& \left(n^{(L)}_1, m^{(L)}_1 \right) 
+  (n_0 +1, m_0+1)  - \left( r^{(L)}_1, t^{(L)}_1\right),     
\end{eqnarray}
%% 
%%%%%%%%%%%%%%%%%%%%%%%%%%%%%%%%%%%%%%%%%%%%%%%%%%%%%%%%%%%%%%%%%%%%%
%%
For the sequence starting from $(4)$ in \eqref{ogromno} we have: 
\begin{eqnarray}
\label{dopcond2}
(n_0, m_0)=(n, m) + \left(n_i^{(RRK_i)}, m_j^{(RRK_j)} \right)
 - \left(r_i^{(RRK_i)}, t_j^{(RRK_j)} \right),   
\end{eqnarray}
for the sequence starting from $(5)$ in \eqref{ogromno}: 
\begin{eqnarray}
\label{dopcond3}
(n, m)=(n_0, m_0) + \left(n_i^{(LLK_i)},
 m_j^{(LLK_j)}\right) - \left(r_i^{(LLK_i)}, t_j^{(LLK_j)}\right).  
\end{eqnarray}
$i$, $j \ge 2$. 
Note that it is assumed that all  
the resulting indices in the compatibility conditions are non-negative. 
For the general complex \eqref{complex}, 
$i \in \Z$. For the bicomplex 
\eqref{Cpq}-\eqref{deltas} we have $p$, $q$ are 
non-negative. 
Since we assume that $0 \le r^{(LLK_i)}_i \le n^{(LLK_i)}_i$, and 
$0 \le t^{(LLK_i)} \le m^{(LLK_i)}_i$, 
 from the compatibility conditions \eqref{condition1}-\eqref{dopcond1} and  
\eqref{dopcond2}-\eqref{dopcond3} 
we see that, 
depending on the signs of  $n - n_0 -\delta_{i,1}$, $m - m_0 -\delta_{i,1}$, and  
only one branch of systems of the form $(2)$-$(2')$ 
exists at each vertex of the three graph 
associated to the double complex \eqref{Cpq}-\eqref{deltas}.   
Recall that according to Theorem 1 of \cite{CM},   
 the definition of the spaces $C^{n, m}$ do not 
depend on the choice of the transversal basis $\U$,  
 still it depends on parameters of foliation $\F$.       
Nevertheless, differential forms $\omega^{n,m}$ 
do depend on holonomy mappings (see subsection{cm})    
 $h_j$, $j \ge 0$, and play the role of extra parameters in the consideration. 
%%

%%%%%%%%%%%%%%%%%%%%%%%%%%%%%%%%%%%%%%%%%%%%%%%%%%%%%%%%%%%%%%%%%%%%%%%%%%%%%%%%%%%%%%%%%%%%%%%%%%
%%
As explained at the end of Section \ref{pseudo}, due to Theorem \eqref{mainpro}, 
 a path $(K_i)$, $i \ge 1$, defining 
 the generators 
$\left\{ \chi, \; \delta^{n_0,m_0}\chi, \; \Phi, \; \delta^{n,m} \Phi, \;  
\alpha_i^{(K_i)}, \;  
\delta^{n(i),m(i)} 
 \alpha_i^{(K_i)} \right\}$, 
$n(i)= n_i^{(K_i)}$, $m(i)= m_i^{(K_j)}$,  
and relations for independent elements of $C^{n, m}$ in \eqref{ogromno},  
form a continual Lie algebra $\mathcal G(\F)$
 with the space of roots provided by the holonomy 
mappings and of parameters  of the foliation $\F$.   
%% 
%%%%%%%%%%%%%%%%%%%%%%%%%%%%%%%%%%%%%%%%%%%%%%%%%%%%%%%%%%%%%%%%%%%%%%%%%%%%%%%%%%%%%%%%%%%
%% 
\subsection{Double cochain complex: Godbillon-Vey type example} 
In this subsection we provide the explicit example for Theorem \eqref{mainpro} 
 in the case of the orthogonality condition \eqref{ortho} applied to the particular case when,  
in the consideration of the previous subsection, $\Phi=\chi$.   
In particular, in differential geometry, 
 the case of a foliation $\F$ of codimension one 
defined by a one-form on a three-dimensional manifold,   
and the formulation of the Godbillon-Vey cohomology class, 
are included in this consideration.   
%%

%%%%%%%%%%%%%%%%%%%%%%%%%%%%%%%%%%%%%%%%%%%%%%%%%%%%%%%%%%%%%%%%%%%%5
%%
We require the orthogonality for $\chi \in C^{n, m}(\F)$ 
within its own bicomplex space, i.e., to satisfy the condition   
\begin{eqnarray}
\label{integrability}
\chi \cdot \delta^{n,m} \chi= 0.   
\end{eqnarray}
%%
%%%%%%%%%%%%%%%%%%%%%%%%%%%%%%%%%%%%%%%%%%%%%%%%%%%%%%%%%%%%%%%%%%%%
%%
First, let us check the condition for \eqref{integrability}. 
Since ${\bf h}_n$ and ${\bf h}_{n+1}$ has to have the same ${\bf h}_n$, 
(according to the derivation of chain condition), thus 
$n+(n+1)-r=0$, but $r=n$, thus, $n+(n+1)-n=0$, i.e., $n=-1$ which is not possible.    
Even if ${\bf h}_n$ and ${\bf h}_{n+1}$ has not to have the same part, 
we obtain $n+(n+1)-r=0$ $r=2n+1$, i.e., $r>n$ which is not also possible. 
%%
%%%%%%%%%%%%%%%%%%%%%%%%%%%%%%%%%%%%%%%%%%%%%%%%%%%%%%%%%%%%%%%%%%%%
%%
Thus, for $\alpha_1^{(R)} \in C^{n',m'}(\F)$ one has   
\begin{eqnarray}
\label{frob}
\delta^{n,m} \chi= \chi \cdot \alpha_1^{(R)}, 
\end{eqnarray}
and
 $n+1=n+n'-r$, $m+1=m+m'-t$, 
 and \eqref{frob} is possible only when $n'=r+ 1$, $m'=t+1$, $0 \le r\le n$, $0 \le t\le n$, and  
$\alpha_1^{(R)} \in C^{r+1,t+1}(\F)$. 
If we require from \eqref{integrability} that for $\alpha_1^{(L)} \in C^{\alpha, \beta} (\F)$,  
$\chi =  \alpha_1^{(L)} \cdot \delta^{n,m} \chi$, 
then $n= \alpha + n +1 -r'$,  $\alpha=r'-1$, i.e., $\alpha$ 
is smaller than the common degree which is not 
possible and thus such $\alpha_1^{(L)}$ does not exist.   
%%
%%%%%%%%%%%%%%%%%%%%%%%%%%%%%%%%%%%%%%%%%%%%%%%%%%%%%%%%%%%%%%%%%%%%%%%%%%%%%%%%%%%%%%%%%%
%% 
Then, as a result of \eqref{ogromno}, we obtain the system of relations: 
%% 
%%%%%%%%%%%%%%%%%%%%%%%%%%%%%%%%%%%%%%%%%%%%%%%%%%%%%%%%%%%%%%%%%%%%%%%%%55
%% 
\begin{eqnarray}
\label{ogromno0}
 0=\chi \cdot \delta^{n,m} \chi,  \qquad (1), \quad 
 0= \left( \delta^{n,m} \chi \right) \cdot \delta^{n,m} \chi \quad  (3),  
\quad \delta^{n,m} \chi= \chi \cdot \alpha^{(R)}_1,  &&  
\end{eqnarray} 
%% 
%%%%%%%%%%%%%%%%%%%%%%%%%%%%%%%%%%%%%%%%%%%%%%%%%%%%%%%%%%%%%%%%%%%%%%%%%%%%%%%%%%%%
%%
 and other relations of the system \eqref{ogromno} at this level  
collapse since further its branches follow from 
 $(3)$ which is trivial.  
%% 

%%%%%%%%%%%%%%%%%%%%%%%%%%%%%%%%%%%%%%%%%%%%%%%%%%%%%%%%%%%%%%%%%%%%%%%%%%%%%%%%%%%%%%%%%%%%%%%%%%%%%%%%%%%%%%%%%%%% 
%%
Let us denote ${\bf h}_n= \left( h_1, \ldots, h_n \right)$, 
 an $n$-tuple of holonomy mappings (see subsection \ref{cm}) 
 (the space of parameters, here given by holonomy data).
%% 
%%%%%%%%%%%%%%%%%%%%%%%%%%%%%%%%%%%%%%%%%%%%%%%%%%%%%%%%%%%%%%%%%%%%%%%%%
%% 
In this example, for forms 
$\left\{\chi({\bf h}_n), \; \delta^{n,m} \chi({\bf h}'_n), \; 
\alpha^{(R)}_1({\bf h}''_{r+1})   
 \right\}$
 define the local part of the following continual Lie algebra by 
identifying the differential forms with generators of $\mathcal G(\F)$ as 
%%
%%%%%%%%%%%%%%%%%%%%%%%%%%%%%%%%%%%%%%%%%%%%%%%%%%%%%%%%%%%%%%%%%%%%%%%%%%%
%% 
\begin{eqnarray}
\label{idento} 
X_{+1}({\bf h}_n)= \chi({\bf h}_n), \;  
X_{-1}({\bf h}'_{n+1})=\delta^{n,m} \chi({\bf h}'_n), \;  
X_0({\bf h}''_{r+1})= \alpha^{(R)}_1({\bf h}''_{r+1}),  
\end{eqnarray}
%%
%%%%%%%%%%%%%%%%%%%%%%%%%%%%%%%%%%%%%%%%%%%%%%%%%%%%%%%%%%%%%%%%%%%%%%%
%%
The next step relation 
%% 
%%%%%%%%%%%%%%%%%%%%%%%%%%%%%%%%%%%%%%%%%%%%%%%%%%%%%%%%%%%
%%
$0=\delta^{n,m} \chi \cdot \alpha^{(R)}_1 
+  \gamma_{n, r+1}\; 
\chi \cdot \delta^{r+1,t+1} \alpha^{(R)}_1$,  
in the system \eqref{ogromno} 
gives rise a higher grading generator identified with  
$\delta^{r+1,t+1} \alpha^{(R)}_1 ({\bf h}''_{r+1})$.  
%%

%%%%%%%%%%%%%%%%%%%%%%%%%%%%%%%%%%%%%%%%%%%%%%%%%%%%%%%%%%%%%%%%%%%%%%%%%%%%%%%%%%%%%%%
%%
According to our derivation of Section \ref{pseudo}, 
though $n$ and $m$ are fixed in the \eqref{idento}, 
the actual tuples ${\bf h}_n$, ${\bf h}'_n$ in $\chi({\bf h}_n)$  
and $\delta^{n,m} \chi({\bf h}'_n)$  
of \eqref{integrability}
 could be different with the same $n$. 
%%

%%%%%%%%%%%%%%%%%%%%%%%%%%%%%%%%%%%%%%%%%%%%%%%%%%%%%%%%%%%%%%%%%%%%%%%%%%%%%%%%%%%%%%%%
%%
Substituting in \eqref{ogromno0} the product $\cdot$ with the graded commutator, we obtain 
 the following commutation relations in the principal grading:  
%%
%%%%%%%%%%%%%%%%%%%%%%%%%%%%%%%%%%%%%%%%%%%%%%%%%%%%%%%%%%%%%%%%%%%%%%
%% 
\begin{eqnarray}
\label{ogromno00}
  && 
\left[X_{+1}({\bf h}_n), X_0({\bf h'}_{n+1})\right]=0,  
 \quad
 \left[X_0({\bf h}_{n+1}), X_0({\bf h}'_{n+1})\right]=0,  
\\   
\nonumber 
%%%%%%%%%%%%%%%%%%%%%%%%%%%%%%%%%%%%%%%%%%%%%%%%%%%%%%%%%%%%%%%%%%% 
%% 
 && \left[X_{+1}({\bf h}_n), X_{-1}(\widetilde {\bf h}_{r+1}) \right]=   
X_0\left( {\bf h}'_{n+1} \right),  
\end{eqnarray}     
where we skip $gr$ in notation of commutators here and what follows. 
%%
%%%%%%%%%%%%%%%%%%%%%%%%%%%%%%%%%%%%%%%%%%%%%%%%%%%%%%%%%%%%%%%%%%%%%%% 
%%
Note also, that the product of two identical differential forms of \cite{CM} vanishes. 
%% 
%%%%%%%%%%%%%%%%%%%%%%%%%%%%%%%%%%%%%%%%%%%%%%%%%%%%%%%%%%%%%%%%%%%%%% 
%%  
Therefore, $[\chi({\bf h}_n), \chi({\bf h}'_n)]=0$, i.e.,  
$[X_{+1}({\bf h}_n), X_{+1}({\bf h}_n)]=0$, and 
$[X_0(\widetilde {\bf h}_{r+1}), X_0(\widetilde {\bf h}'_{r+1})]=0$.    
%% 
%%%%%%%%%%%%%%%%%%%%%%%%%%%%%%%%%%%%%%%%%%%%%%%%%%%%%%%%%%%%%%%%%%%%%%%%%
%%
In this example we use a non-standard grading. 
%%
%%%%%%%%%%%%%%%%%%%%%%%%%%%%%%%%%%%%%%%%%%%%%%%%%%%%%%%%%%%%%%%%%%
%% 
Note that an element of $C^{n,m}(\F)$ is an $n$-form
 $\omega(h_1, \ldots, h_{n})$ depending on $n$ holonomy maps. 
%% 
%%%%%%%%%%%%%%%%%%%%%%%%%%%%%%%%%%%%%%%%%%%%%%%%%%%%%%%%%%%%%%%%%%%%%%%%% 
%%  
Thus, the space of continual roots 
is provided by the space of holonomy embeddings (see subsection \ref{cm}).   
%% 
%%%%%%%%%%%%%%%%%%%%%%%%%%%%%%%%%%%%%%%%%%%%%%%%%%%%%%%%%%%%%%%%%%%%%%%%%%%%%%%%%%%%%%
%% 
From \eqref{ogromno00} we extract directly the kernels %(using the linearity)   
%%
%%%%%%%%%%%%%%%%%%%%%%%%%%%%%%%%%%%%%%%%%%%%%%%%%%%%%%%%%%%%%%%%%%%%%%%%%%%%%%%%%%%%
%% 
\begin{eqnarray*} 
 K_{+1, 0}\left({\bf h}_n, {\bf h}'_{r+1} \right)=0, 
\quad  K_{0, 0 }\left({\bf h}_n, {\bf h}'_{r+1} \right)=0, \quad 
 K_{+1, -1}\left({\bf h}_n, \widetilde {\bf h}_{r+1} \right) 
=  {\bf h}'_{n+1},  
\end{eqnarray*}  
%%
%%%%%%%%%%%%%%%%%%%%%%%%%%%%%%%%%%%%%%%%%%%%%%%%%%%%%%%%%%%%%%%%%%%%%%
%% 
Therefore,   
$K_{+1, +1}({\bf h}_n, {\bf h}'_n)=0$ and 
%%
%%%%%%%%%%%%%%%%%%%%%%%%%%%%%%%%%%%%%%%%%%%%%%%%%%%%%%%%%%%%%%%%%%%%%%%%%%%%%%%%  
%%
$K_{0, 0}({\bf h}''_{r+1}, {\bf h}''_{r+1})=0$.  
%%  
%%%%%%%%%%%%%%%%%%%%%%%%%%%%%%%%%%%%%%%%%%%%%%%%%%%%%%%%%%%%%%%%%%%%%%%%%%%%%%%%%%%%%%%%%%
%%
In the case of a codimension one foliation defined by a one-form,  
$\chi$, $\alpha^{(R)}_1$, $\chi \in C^{1,m}(\F)$.    
Recall \cite{Ghys} that the Godbillon-Vey cohomology class is given by 
$\left[ \alpha^{(R)}_1 \wedge \delta^{1,m} \alpha^{(R)}_1\right]$.
The construction above clarifies the Lie-algebraic meaning of this cohomology class. 
%%

%%%%%%%%%%%%%%%%%%%%%%%%%%%%%%%%%%%%%%%%%%%%%%%%%%%%%%%%%%%%%%%%%%%%%%%%%%%%%%%%%%%%%%%%%%%%%%
%%
We then check the graded the Jacobi identity  
 for generators \eqref{idento}.    
%%
%%%%%%%%%%%%%%%%%%%%%%%%%%%%%%%%%%%%%%%%%%%%%%%%%%%%%%%%%%%%%%%%%%%%%%%%%%%%%%%%%%%%
%%
For instance, let us do that  
for the generators 
$X_{+1}({\bf h}_n)$, 
$X_{-1} (\widetilde {\bf h}_{r+1}) $, and $X_0\left({\bf h}'_{n+1}\right)$.  
%%  
%%%%%%%%%%%%%%%%%%%%%%%%%%%%%%%%%%%%%%%%%%%%%%%%%%%%%%%%%%%%%%%%%%%5
%%
From \eqref{ogromno00} we obtain 
\[
\left[X_{-1}(\widetilde {\bf h}_{n+1}), X_0({\bf h}'_{r+1}) \right]  
= -\gamma_{n,r+1} 
\left[X_{+1}(\widetilde {\bf h}_n), Y_0({\bf h}'_{r+2}) \right]. 
\]
%% 
%%%%%%%%%%%%%%%%%%%%%%%%%%%%%%%%%%%%%%%%%%%%%%%%%%%%%%%%%%%%%%%%%%%%
%% 
 Usiing the formula for the double graded commutator, we see that 
%% 
%%%%%%%%%%%%%%%%%%%%%%%%%%%%%%%%%%%%%%%%%%%%%%%%%%%%%%%%%%%%%%%%%%%%
%% 
\begin{eqnarray*}
&&\left[ X_{+1}({\bf h}_n), 
\left[X_{+1}({\bf h}_n), Y_{0}({\bf h}'_{r+2}) \right]\right] 
\nn
&&
= \left[ \left[ X_{+1}({\bf h}_n),  
X_{+1}({\bf h}_n)\right], Y_{0}({\bf h}'_{r+2}) \right] 
+ (-1)^{n n}  \left[ X_{+1}({\bf h}_n),  
\left[X_{+1}({\bf h}_n), Y_{0}({\bf h}'_{r+2})\right]\right],  
\end{eqnarray*}
%%
%%%%%%%%%%%%%%%%%%%%%%%%%%%%%%%%%%%%%%%%%%%%%%%%%%%%%%%%%55
%%
 when $n^2=2k+1$, $k\in \mathbb Z$, 
$\left[  X_{+1}({\bf h}_n), \left[ X_{-1}(\widetilde {\bf h}_{r+1}), 
X_{0}\left( {\bf h}'_{n+1} \right) \right] \right]=0$.  
%%  
%%%%%%%%%%%%%%%%%%%%%%%%%%%%%%%%%%%%%%%%%%%%%%%%%%%%%%%%%%%%%%%%%%%%%%%%
%% 
Then it follows 
%% 
%%%%%%%%%%%%%%%%%%%%%%%%%%%%%%%%%%%%%%%%%%%%%%%%%%%%%%%%%%%%%%%%%%%%%%%%%
%% 
\begin{eqnarray*}
&&\gamma_{n+1,r+1}\; 
\left[ X_{0}({\bf h'}_{n+1}), \left[X_{+1}({\bf h}_n), 
X_{-1}(\widetilde {\bf h}_{r+1}) \right] \right] 
\nn
&&
+\gamma_{n,n+1}\; \left[  X_{+1}({\bf h}_n), \left[  
X_{-1}(\widetilde {\bf h}_{r+1}),  X_{0}({\bf h'}_{n+1}) \right] \right] 
\nn
&&+\gamma_{r+1, n}\; \left[X_{-1}(\widetilde {\bf h}_{r+1}),  \left[ 
X_{0}({\bf h'}_{n+1}), X_{+1}({\bf h}_n) \right] \right] 
%%
%%%%%%%%%%%%%%%%%%%%%%%%%%%%%%%%%%%%%%%%%%%%%%%%%%%%%%%%%%%% 
\nn
&&=
\gamma_{n+1,r+1}  
\left[ X_{0}({\bf h'}_{n+1}), X_{0}({\bf h'}_{n+1}) \right]  
+ \gamma_{n,n+1} \gamma_{n, r+1} \left[  X_{+1}({\bf h}_n),    
\left[ X_{+1}({\bf h}_n), Y_{0}({\bf h}'_{r+2}) \right] \right]  
\nn
&&+\gamma_{r+1, n}\; \left[X_{-1}(\widetilde {\bf h}_{r+1}), 0  \right]=0, 
\end{eqnarray*}
and the Jacobi identity is satisfied. 
%%  
%%%%%%%%%%%%%%%%%%%%%%%%%%%%%%%%%%%%%%%%%%%%%%%%%%%%%%%%%%%%%%%%%%%%%%5
%%
 Thus, the generators defined in \eqref{idento} 
satisfy the Jacobi identity, 
confirming the consistency of the continual Lie algebra 
$\mathcal{G}(\mathcal{F})$ associated with the chain complex.
%% 
%%%%%%%%%%%%%%%%%%%%%%%%%%%%%%%%%%%%%%%%%%%%%%%%%%%%%%%%%%%%%%%%%%%%%%%%%%%%%%%%%%%%%%%%%%
%%%%%%%%%%%%%%%%%%%%%%%%%%%%%%%%%%%%%%%%%%%%%%%%%%%%%%%%%%%%%%%%%%%%%%%%%%%%%%%%%%%%%%%%%%
%% 
\section{Conclusions}
In conclusion, we would like to mention a few directions of development and 
further applications of the material presented in this paper. 
We propose a way to associate a continual Lie algebra to a chain complex. 
Thus the properties, in particular,  the Jacobi identity, kernels, 
and relations of resulting 
continual Lie algebras   
depend on the kind of chain complex spaces as well as on the set of their parameters.  
%% 
%%%%%%%%%%%%%%%%%%%%%%%%%%%%%%%%%%%%%%%%%%%%%%%%%%%%%%%%%%%%%%%%%%%%%%%%%%%%%%%%%%%%%%%%%%
%%
 One can think of introducing various types of products suitable for the construction of 
systems of relations 
more complicated than \eqref{ogromno}.   
In our particular case (Section \ref{pseudo}), 
in order to make connection with continual Lie algebras, 
 we have chosen the commutator \eqref{product0} 
(with respect to the original product 
defined on bicomplex spaces) as the simplest natural product.  
One could think of other possibilities which would be  
 coherent with the orthogonality condition \eqref{ortho}. 
%% 

%%%%%%%%%%%%%%%%%%%%%%%%%%%%%%%%%%%%%%%%%%%%%%%%%%%%%%%%%%%%%%%%%%%%%%%%
%% 
In our exposition, the standard form of chain complexes was involved.  
Nevertheless, one can consider more complicated setups, 
 in particular, chain complexes where differentials
 act in non-trivial ways with respect to 
indices of spaces (c.f. examples in \cite{Huang}).   
%%
%%%%%%%%%%%%%%%%%%%%%%%%%%%%%%%%%%%%%%%%%%%%%%%%%%%%%%%%%%%%%%%%%%%%%%%
%%
That would lead to alternative forms of systems of relations 
as well as of compatibility conditions. 
What could be especially interesting, 
is to treat multiple chain complexes containing combinations of a few chain-cochains.  

%%%%%%%%%%%%%%%%%%%%%%%%%%%%%%%%%%%%%%%%%%%%%%%%%%%%%%%%%%%%%%%%%%%%%%%%%%
%%
The example of the \v{C}hech-de Rham cohomology that we study in Section \ref{example} 
 comes from the differential geometry of foliations. 
In classics, the orthogonality condition applied 
to elements and their differentials of one 
particular bicomplex space, 
boils down to the integrability condition, and leads to the Frobenius theorem. 
Then it delivers the Godbillon-Vey cohomology  
class whose geometric meaning is not yet completely studied 
\cite{galaev}.   
%%
%%%%%%%%%%%%%%%%%%%%%%%%%%%%%%%%%%%%%%%%%%%%%%%%%%%%%%%%%%%%%%
%%
As for further applications in differential geometry,
starting from the orthogonality condition,  
it would be interesting to find other cohomology 
invariants related to the \v{C}ech-de Rham 
bicomplex for foliations,
 and to understand their geometric meaning \cite{Zu100, Zu, zufo}.
%%
%%%%%%%%%%%%%%%%%%%%%%%%%%%%%%%%%%%%%%%%%%%%%%%%%%%%%%%%%%%%%%%%%%%%     
%% 
The constructions considered in this paper can be also used for 
the cohomology theory of smooth manifolds \cite{Zu4, Zu1, Zu3, Zu2},
 in particular, in 
 various approaches to the construction of cohomology classes \cite{Zu6, Zu5}  
(cf. \cite{Losik}).  
Wide applications are awaiting for new examples of 
continual Lie algebras in the field of 
integrable models 
\cite{lezsav, razsav}.       
%%%
The cases of non-commutative fields used to 
define continual Lie algebras would also be useful for in  
non-commutative geometry.  
%%

%%%%%%%%%%%%%%%%%%%%%%%%%%%%%%%%%%%%%%%%%%%%%%%%%%%%%%%%%
%%
Concerning applications in Physics,
 the derivation of continual Lie algebras
given in this paper, 
will be also useful in mathematical physics describing condensed matter theory. 
%%
%%%%%%%%%%%%%%%%%%%%%%%%%%%%%%%%%%%%%%%%%%%%%%%%%%%%%%%%%
%% 
Let us mention 
non-perturbative, and topological defects dominates dynamics \cite{z8},
fermionic superfluids \cite{z7}, 
Wigner-Weyl calculus \cite{z1, z2}, 
the theory of topological invariants and relations  
relation between solid state systems and high energy physics \cite{z4, z5, z6},  
chiral separation effect \cite{z3}, 
and non-renormalization by interactions of integer quantum Hall effect \cite{z9}.  

%%
%%%%%%%%%%%%%%%%%%%%%%%%%%%%%%%%%%%%%%%%%%%%%%%%%%%%%%%%%%%%%%%%%%%%%%%%%%%%%%%%%%%%%%%%%%
%%%%%%%%%%%%%%%%%%%%%%%%%%%%%%%%%%%%%%%%%%%%%%%%%%%%%%%%%%%%%%%%%%%%%%%%%%%%%%%%%%%%%%%%%%
%%
\section*{Acknowledgments}
The author would like to thank an unknown referee for valuable corrections and comments. 
The author is supported by the Institute of Mathematics, 
Academy of Sciences of the Czech Republic (RVO 67985840), and  
Scientific Collaboration Grant with Researchers from the Czech Republic,   
Israel Science Foundation.
%%
%%%%%%%%%%%%%%%%%%%%%%%%%%%%%%%%%%%%%%%%%%%%%%%%%%%%%%%%%%%%%%%%%%%%%%%%%%%%%%%%%%%%%%%%%%
%%%%%%%%%%%%%%%%%%%%%%%%%%%%%%%%%%%%%%%%%%%%%%%%%%%%%%%%%%%%%%%%%%%%%%%%%%%%%%%%%%%%%%%%%%
%%  

%% 
\end{document}